\documentclass[hidelinks,onefignum,onetabnum]{siamart251216}

\usepackage{lipsum}
\usepackage{amsfonts}
\usepackage{graphicx}
\usepackage{algorithmic}
\usepackage{amssymb}
\usepackage{array}
\usepackage{booktabs}
\usepackage{multirow}

\ifpdf
  \DeclareGraphicsExtensions{.pdf,.png,.jpg}
\else
  \DeclareGraphicsExtensions{.eps}
\fi

\newsiamremark{remark}{Remark}
\newsiamremark{hypothesis}{Hypothesis}
\crefname{hypothesis}{Hypothesis}{Hypotheses}
\newsiamthm{claim}{Claim}
\newsiamremark{fact}{Fact}
\crefname{fact}{Fact}{Facts}

\headers{Numerical Approximation of MTEP in Dispersive Media}{C. Wang}

\title{Numerical Approximation of Maxwell Transmission Eigenvalue Problem in Dispersive Media via Holomorphic Operator Function\thanks{
\funding{The work was supported by National Natural Science Foundation of China grant No. 12501558.}}}

\author{Chao Wang\thanks{MSU-BIT-SMBU Joint Research Center of Applied Mathematics, Shenzhen MSU-BIT University, Shenzhen 518001, China. 
  (\email{wangchao@smbu.edu.cn}).}
}

\usepackage{amsopn}

\def\cT{\mathscr{T}}

\def\O{\Omega}
\def\p{\partial}

\def\bpsi{\boldsymbol{\psi}}
\def\bphi{\boldsymbol{\phi}}
\def\bxi{\boldsymbol{\xi}}
\def\bmu{\boldsymbol{\mu}}

\def\bu{\boldsymbol{u}}

\def\bv{\boldsymbol{v}}
\def\bw{\boldsymbol{w}}

\def\bn{\boldsymbol{n}}

\def\be{\boldsymbol{e}}
\def\bof{\boldsymbol{f}}
\def\curl{\nabla\times}
\def\grad{\nabla}
\def\div{\nabla\cdot}

\def\cT{\mathcal{T}}

\def\bn{\boldsymbol{n}}

\def\norm#1{|\hskip -1pt|#1|\hskip -1pt|}

\def\bn{\boldsymbol{n}}
\def\bg{\boldsymbol{g}}

\def\bA{\boldsymbol{A}}
\def\bN{\boldsymbol{N}}
\def\bx{\boldsymbol{x}}

\ifpdf
\hypersetup{
  pdftitle={Numerical Approximation of Maxwell Transmission Eigenvalue Problem in Dispersive Media via Holomorphic Operator Function},
  pdfauthor={Chao Wang}
}
\fi

\begin{document}

\maketitle

\begin{abstract}
In this paper, we develop and analyze a numerical method for the Maxwell's transmission eigenvalue problem in dispersive media. The problem is formulated as the eigenvalue problem of a holomorphic Fredholm operator function of index zero. For discretization, we employ the N\'ed\'elec elements. Convergence of the scheme is established via the abstract approximation theory of Karma \cite{karma1996approximationI, karma1996approximationII}. To compute the resulting nonlinear matrix eigenvalue problems, we employ a parallel multi-step contour integral method. Numerical experiments are presented to demonstrate the viability and accuracy of the proposed approach.
\end{abstract}



\section{Introduction}\label{sec:int}
The transmission eigenvalue problem (TEP) plays an important role in the inverse scattering theory of inhomogeneous media. It appears in various qualitative methods and can be used to estimate material properties using scattering data~\cite{CCDMS:2010:InvElemgSctPlmAstMd, CCH:2010:ODtmDrcTsmEgivFFldDt, SJ:2010:EstTsmEgivIdRftCchyDt, SJ:2012:EgnvMthdUMlpFqcDtInvSctPlm, CCM:2011:LnSmpMthdInvElemgSct}. The TEP and its variants have attracted considerable attention, leading to substantial advancements in their theoretical analysis and practical applications \cite{CK:1998:InvAcuElemgSctThy, dhia2011use, CH:2012:TsmEglvIvsSctThy, KA:2009:OnExtTsmEigl, PS:2008:TE, CGH:2010:ExtIftDscStTsmEgnl}. There has also been an increasing interest within the scientific computing community in developing new numerical methods for TEPs~\cite{CMS:2010:AnyCmpMthdTsmEgnv, HLSX:2011:ACpBFItTsmPlmAcst, CPS:2007:IntTsmPlm, mora2018virtual, sun2016finite, xi2020multi, xie2017multilevel, YGB:2016:NCFEMTEP, YBLH:2016:MMHTE, KA:2013:AnumericalMCITE}.

In this paper, we consider the computation of the Maxwell's transmission eigenvalue problem (MTEP) in dispersive media. Let $\O\subset \mathbb{R}^3$ be a bounded simply connected domain with Lipschitz boundary $\partial \Omega$ occupied by a dispersive and anisotropic medium. Let $\bA(\bx,k), \bN(\bx,k)\in L^{\infty}(\O,\mathbb{C}^{3\times 3})$ be Hermitian matrix-valued functions satisfying $A_*|\bxi|^2<\overline{\bxi}\cdot\bA\bxi<A^*|\bxi|^2$ and $N_*|\bxi|^2<\overline{\bxi}\cdot\bN\bxi<N^*|\bxi|^2$ for all $\bxi\in\mathbb{C}^3$ a.e. in $\O$ and $k\in\mathbb{R}$, where $A_*$, $A^*$, $N_*$ and $N^*$ are positive constants independent of $k$. We assume that one of the following cases holds uniformly with respect to $k$: 1) $A_*>1$ and $N_*>1$, 2) $A_*>1$ and $N^*<1$, 3) $A^*<1$ and $N^*<1$, 4) $A^*<1$ and $N_*>1$. In addition, $\bA(\cdot,k)$ and $\bN(\cdot,k)$ are assumed to be holomorphic with respect to $k\in\mathcal D\subset\mathbb{C}$. The MTEP is to find $k^2\in\mathbb{C}\backslash\{0\}$ and non-trivial functions $\bu\in H(\text{curl},\O)$ and $\bv\in H(\text{curl},\O)$ such that
\begin{equation}\label{MTEPA}
\left\{
\begin{aligned}
\curl\left( \bA(\bx,k)\curl\bu \right) - k^2\bN(\bx,k)\bu &= \boldsymbol{0}\quad\text{in}\,\O, \\
\curl\curl\bv - k^2\bv &= \boldsymbol{0} \quad\text{in}\,\O,\\
(\bu-\bv)\times\bn &= \boldsymbol{0} \quad\text{on}\,\p\O, \\
(\bA(\bx,k)\curl\bu)\times\bn -(\curl\bv)\times\bn &=\boldsymbol{0} \quad\text{on}\,\p\O,
\end{aligned}
\right.
\end{equation}
where $\bn$ is the unit outward normal to $\p\O$. Here, $\bA(\bx,k)$ and $\bN(\bx,k)$ correspond to the permeability and permittivity of the medium, respectively. The dependence of $\bA(\bx,k)$ and $\bN(\bx,k)$ on the wave number $k$ models the dispersive response of the medium. When $\bA$ and $\bN$ are independent of $k$, problem~\eqref{MTEPA} reduces to the standard anisotropic problem. The MTEP consists of a system of second-order partial differential equations involving vector-valued fields, which is nonlinear and non-selfadjoint. In the dispersive case, the dependence of the material parameters $\bA(\bx,k)$ and $\bN(\bx,k)$ on the spectral parameter introduces an additional nonlinearity into the problem. Numerical discretization for \eqref{MTEPA} therefore leads naturally to nonlinear, non-Hermitian matrix eigenvalue problems. Computing eigenvalues in such cases poses significant challenges, particularly for large-scale matrices and when a priori information on the spectral distribution is unavailable. 

Numerical methods for the Helmholtz transmission eigenvalue problem have been proposed by many researchers (see, e.g., \cite{CMS:2010:AnyCmpMthdTsmEgnv, JLS:2020:CIETC, SJ:2011:ItvMthdTsmEgnv,JSX:2014:MGMthdHmtzTsmEgnvPlm, JS:2013:MLMTEAM, CMS:2014:EAFEATE,YGB:2016:NCFEMTEP,GJSX:2016:C0IPTEP,YBLH:2016:MMHTE,BHY:2021:LPFESETEP,XJG:2018:C0IPTEEW}). In contrast, far fewer results are available for the Maxwell case. Monk and Sun \cite{MS:FEMMxwTsmEgnv} introduced the curl-conforming method and the mixed finite-element method. Sun and Xu~\cite{SunXu:2013:CMTEAIMP} proposed an iterative method based on the fourth-order formulation. An and Zhang~\cite{AZ:2018:EfcSpctGlkAppxmEEMxwTsmEgnvSphGeo} investigated the spectral method for the MTEP on a spherical region (see also \cite{An2024SIAM}). Huang et al.~\cite{HHL:2015:ArobustNAMTEP} developed an eigenvalue solver for computing a few of the smallest positive transmission eigenvalues. Recently, Han and Zhang~\cite{HZ:2022:Hcurl2CEMTEPFA} introduced a fixed-point method using $H(\text{curl}^2)$-conforming elements. Liu et al.~\cite{liu2025mixed} developed a new mixed finite element scheme. In addition, Han~\cite{han2024discontinuous} proposed a DG method. In \cite{WCS:2024:EEMFEMMTEP}, we proved error estimates for the iterative method proposed in \cite{SunXu:2013:CMTEAIMP}. Cakoni et al. proposed a computational method based on boundary integral equations \cite{Cakoni2018}. All of these studies concern frequency-independent material parameters, whereas numerical approximation of the MTEP with dispersive coefficients remains unexplored.

The error analysis for the MTEP is more challenging than that for the Helmholtz TEP. Existing numerical results have been predominantly reported for isotropic media, i.e., $\bA=\boldsymbol{I}$ and $\bN=\bN(\bx)$. In such cases, the MTEP can be reformulated as a nonlinear fourth-order problem. Thus one can derive the regularity of the solution space and the compactness property of the solution operator. Consequently,  convergence analysis can be carried out using classical spectral perturbation theory for linear operators. In contrast, for dispersive media when $\bA=\bA(\bx,k)$ and $\bN=\bN(\bx,k)$, the fourth-order reformulation is no longer available. To the best of our knowledge, numerical methods for the dispersive MTEP remain limited. Recent contributions are only for the anisotropic cases which are given by \cite{han2025finite, HSZ:ModifyMTEP:2026}. Moreover, in dispersive media the associated operator depends nonlinearly on $k$. Consequently, the classical approximation theory based on a fixed compact solution operator cannot be applied directly.

Motivated by~\cite{gong2022finite, xi2022holomorphic}, we propose and analyse a finite element method  for the dispersive MTEP based on a holomorphic operator function framework. The holomorphic dependence of $\bA(\bx,k)$ and $\bN(\bx,k)$ on the wave number makes this framework particularly suitable for treating frequency-dependent material parameters. Specifically, we reformulate \eqref{MTEPA} as the eigenvalue problem of a holomorphic Fredholm operator function. The  N\'ed\'elec edge elements are used to discretize the resulting operator function. Using the parallel multi-step contour integral method, we compute eigenvalues inside a prescribed region of the complex plane. Furthermore, combining standard finite element approximation theory with the abstract eigen-approximation theory for holomorphic Fredholm operator functions, we establish error estimates for the computed eigenvalues and the associated generalized eigenvectors.

The rest of the paper is organized as follows. In \cref{sec:prelimi}, we introduce some preliminaries of holomorphic Fredholm operators and the abstract approximation theory for the eigenvalue problem. In \cref{sec:optTEP}, we reformulate the MTEP as an operator eigenvalue problem and show that the operator is a holomorphic Fredholm operator of index zero. The \cref{sec:FEP} contains the finite element discretization of the operator functions and the convergence of the eigenvalue problem. In \cref{sec:pmCIMandNumExp}, a parallel multi-step contour integral method is designed to compute the eigenvalues and numerical results are presented. Finally, we make some conclusions and discuss future work in \cref{sec:con}.

\section{Preliminaries}\label{sec:prelimi}                  

We present the approximation theory for eigenvalues of holomorphic Fredholm operator functions~\cite{beyn2014finding, karma1996approximationI, karma1996approximationII}. Let $X$ and $Y$ be complex Banach spaces. Denote by $\mathcal{L}(X,Y)$ the space of bounded linear operators from $X$ to $Y$. For $A\in\mathcal{L}(X,Y)$, let $\mathcal{N}(A)=\{ x\in X: Ax=0 \}$ and $\mathcal{R}(A)=\{ y\in Y: y=Ax,\,x\in X \}$ be the null space and the range of $A$, respectively. Let $\mathcal{D}\subset\mathbb{C}$ be a bounded simply connected region and $F:\mathcal{D}\rightarrow\mathcal{L}(X,Y)$ be a holomorphic operator function. 

\begin{definition}
A bounded linear operator $F: X\rightarrow Y$ is said to be Fredholm if 
\begin{itemize}
\item[1.] the subspace $\mathcal{R}(F)$, range of $F$, is closed in $Y$;
\item[2.] the subspace $\mathcal{N}(F)$, null space of $F$, and $Y/\mathcal{R}(F)$ are finite-dimensional.
\end{itemize}
The index of $F$ is the integer defined by 
$
ind(F)=\text{dim}\mathcal{N}(F)-\text{dim}(Y/\mathcal{R}(F)).
$
\end{definition}

In the rest of the paper, we assume that $F(\eta)$ is a holomorphic operator function and for each $\eta\in\O$, $F(\eta)$ is a Fredholm operator of index $0$. 

\begin{definition}
A complex number $\lambda\in\mathcal{D}$ is called an eigenvalue of $F$ if there exists a nontrivial $x\in X$ such that $F(\lambda)x=0$. The element $x$ is called an eigenelement associated with $\lambda$. 
\end{definition}

The resolvent set $\rho(F)$ and the spectrum $\sigma(F)$ of $F(\cdot)$ are defined as
\begin{align*}
\rho(F)=\{ \eta\in\mathcal{D}: F(\eta)^{-1}\,\,\text{exists}\,\,\text{and}\,\, \text{is}\,\, \text{bounded} \}\quad\text{and}\quad \sigma(F)=\mathcal{D}\backslash\rho(F),
\end{align*}
respectively. If $\rho(F)\neq \emptyset$, since $F(\eta)$ is holomorphic, the spectrum $\sigma(F)$ has no cluster points in $\mathcal{D}$ and every $\lambda\in\sigma(F)$ is an eigenvalue for $F$. Furthermore, the operator function $F^{-1}(\cdot)$ is meromorphic~\cite{karma1996approximationI}. The dimension of $\mathcal{N}(F(\lambda))$ is called the geometric multiplicity of an eigenvalue $\lambda$. 

\begin{definition}
An ordered sequence of elements $x_0,\, x_1,\cdots, x_k$ in $X$ is called a Jordan chain of $F$ at an eigenvalue $\lambda$ if
\begin{align*}
F(\lambda)x_j + \frac{1}{1!}F^{(1)}(\lambda)x_{j-1}+\cdots +\frac{1}{j!}F^{(j)}(\lambda)x_0=0,\quad j=0,1,\cdots,k,
\end{align*}
where $F^{(j)}$ denotes the $j$th derivative.
\end{definition}

The length of any Jordan chain of an eigenvalue is finite. Denote by $m(F,\lambda,x_0)$ the maximal length of a Jordan chain initiated by $x_0$. The maximal length of all Jordan chains of the eigenvalue $\lambda$ is denoted by $\kappa(F,\lambda)$. Elements of any Jordan chain of an eigenvalue $\lambda$ are called generalized eigenelements of $\lambda$.

\begin{definition}
The linear span of all generalized eigenelements of an eigenvalue $\lambda$, denoted by $G(\lambda)$, is called the generalized eigenspace of $\lambda$.
\end{definition}

%

Let $X_n, Y_n$, $n\in\mathbb{N}$ by general Banach spaces. Denote by $\Phi_0(X,Y)$ the sets of all Fredholm operators of index zero in $\mathcal{L}(X,Y)$. Consider a sequence of holomorphic Fredholm operator functions such that
\begin{align*}
F_n:\mathcal D\rightarrow \Phi_0(X_n,Y_n),\quad n\in \mathbb{N}.
\end{align*}

Assume that $\rho(F)\neq\emptyset$ and the following properties hold
\begin{itemize}
\item[(a1)] There exist bounded linear mappings $p_n\in\mathcal{L}(X, X_n)$ and $q_n\in \mathcal{L}(Y, Y_n)$ be such that 
\begin{align}
\lim_{n\rightarrow\infty}\norm{p_nx}_{X_n}=\norm{x}_X,\,\,\, x\in X,\quad\text{and}\quad \lim_{n\rightarrow\infty}\norm{q_ny}_{Y_n}=\norm{y}_Y,\,\,\, y\in Y.
\end{align}
\item[(a2)] The sequence $\{F_n(\cdot)\}$ is equibounded on $\mathcal{D}$, i.e., $\norm{F_n(\eta)}\leq C,\,\eta\in\mathcal{D}$;
\item[(a3)] $\{F_n(\eta)\}_{n\in \mathbb{N}}$ converges to $F(\eta)$ for every $\eta\in\mathcal{D}$ and $x\in X$, i.e.,
\begin{align*}
\norm{[F_n(\eta)p_n-q_nF(\eta)]x}_{Y_n}\rightarrow 0;
\end{align*}
\item[(a4)] For any subsequence $x_n\in X_n$, $n\in N^\prime$ with $\norm{x_n}_{X_n}$, $n\in N^\prime$ bounded and $\lim_{N^\prime\ni n\rightarrow\infty}\norm{F_n(\eta)x_n-q_n y}_{Y_n}=0$ for some $y\in Y$, there exists a subsequence $N^{\prime\prime}\subset N^\prime$ such that 
\begin{align*}
\lim_{N^{\prime\prime}\ni n\rightarrow\infty}\norm{x_n-p_nx}_{Y_n}=0.
\end{align*}
\end{itemize}

The following abstract approximation results hold true~\cite{karma1996approximationI, karma1996approximationII, gohberg2009holomorphic}.
\begin{theorem}\label{thm:abstractapprholoFredmOpt}
Assume that (a1)-(a4) hold. For any $\lambda\in\sigma(F)$ there exists $n_0\in \mathbb{N}$ and a sequence $\lambda_n\in\sigma(F_n)$, $n\geq n_0$, such that $\lambda_n\rightarrow\lambda$ as $n\rightarrow\infty$. For any sequence $\lambda_n\in\sigma(F_n)$ with this convergence property, and associated eigenvectors $v_n^0\in\mathcal{N}(F_n(\lambda_n)),\, \norm{v_n^0}_{X_n}=1$,
one has that 
\begin{align*}
|\lambda_n-\lambda| &\leq C\epsilon_n^{1/\kappa}, \\
\inf_{v\in\mathcal{N}(F(\lambda))}\norm{v_n^0-p_nv}_{X_n} &\leq C\epsilon_n^{1/\kappa},
\end{align*}
where $\epsilon_n=\max_{|\eta-\lambda|\leq\delta}\max_{v\in G(\lambda)}\norm{F_n(\eta)p_nv-q_nF(\eta)v}_{Y_n}$, for sufficiently small $\delta>0$.
\end{theorem}

\section{Transmission Eigenvalue Problem}\label{sec:optTEP}                  

In this section, we formulate the problem \eqref{MTEPA} as the eigenvalue problem of a holomorphic operator function. For simplicity, we set $\eta=k^2$ and consider the separable dispersive coefficients $\bA(\bx,\eta)=q_A(\eta)\bA_0(\bx)$ and $\bN(\bx,\eta)=q_N(\eta)\bN_0(\bx)$ where $\bA_0(\bx),\,\bN_0(\bx)\in L^\infty(\O,\mathbb{C}^{3\times 3})$ are Hermitian as well as uniformly positive definite. Moreover, $q_A(\eta)$, $q_N(\eta)\neq0$ in $\mathbb{C}$ and $q_A(\eta),\, q_N(\eta)>0$ for positive real values of $\eta$.
Define
\begin{align*}
L_t^2(\p\O) &=\{ \bg\in [L^2(\p\O)]^3:\,\bg\cdot\bn=0\quad\text{a.e.\,\,on}\,\,\p\O  \}, \\
X &=\{ \bu\in H(\text{curl},\O):\,\bu\times\bn\in [L^2(\p\O)]^3 \}, \\
S &=\{ p\in H^1(\O): p|_{\p\O} \,\,\text{is\,\,a\,\,constant\,\,and}\,\,\int_{\O} p\,\mathrm{d}x=0 \}.
\end{align*}
The space $X$ is equipped with the norm $\norm{\bu}^2_X=\norm{\bu}^2_{H(\text{curl},\O)}+\norm{\bu\times\bn}^2_{L^2(\p\O)}.$ Besides, the following weighted Helmholtz decomposition holds for $X$ \cite{Monk:2003:FEMMWE}
\begin{align}\label{eq:decompositionforX}
X=X_{N_0}\oplus\grad S,
\end{align}
where $X_{N_0} =\{ \bu\in X: (\bN_0\bu,\grad\xi)=0,\,\,\forall\,\xi\in S \}$. Thus, $X_{N_0}$ is a closed subspace of $X$ endowed with the norm inherited from $X$.

We consider an auxiliary Maxwell equation. Given $\bg\in L_t^2(\p\O)$, find $\bu$ such that
\begin{equation}\label{def:MaxwellwithAN}
\left\{
\begin{aligned}
\curl\left( q_A(\eta)\bA_0\curl\bu \right) - p(\eta)\bN_0\bu &= \boldsymbol{0}\quad\text{in}\,\,\O, \\
\left(q_A(\eta)\bA_0\curl\bu\right)\times\bn - i(\bn\times\bu)\times\bn &=\bg \quad\text{on}\,\,\p\O,
\end{aligned}
\right.
\end{equation}
where $p(\eta)=\eta q_N(\eta)$. Let $\displaystyle(\bu,\bv):=\int_{\O}\bu\cdot\overline\bv\,\mathrm{d}x$ and $\displaystyle\left\langle\bu,\bv\right\rangle:=\int_{\p\O}\bu\cdot\overline{\bv}\,\mathrm{d}s$. The weak formulation of \eqref{def:MaxwellwithAN} is to find $\bu\in X$ such that
\begin{align}\label{eq:curlcurlequ}
a(\bu,\bv)=\ell(\bv) \quad \forall\, \bv\in X,
\end{align}
where $a(\bu,\bv)$ and $\ell(\bv)$ are defined, respectively, as
\begin{align}
a(\bu,\bv) &=q_A(\eta)\left(\bA_0\curl\bu,\curl\bv\right)-p(\eta)\left(\bN_0\bu,\bv\right)-i\left\langle\bu\times\bn,\bv\times\bn\right\rangle, \label{eq:defa} \\[8pt]
\ell(\bv) &= \left\langle\bg,\gamma_T\bv\right\rangle\quad\text{with}\quad\gamma_T\bv:=(\bn\times\bv)\times\bn.
\end{align}
Due to \eqref{eq:decompositionforX}, we can decompose the solution $\bu=\bu_0+\grad p$ with $\bu_0\in X_{N_0}$ and $p\in S$. Then \eqref{eq:curlcurlequ} is equivalent to the problem of finding $\bu\in X_{N_0}$ such that
\begin{align}\label{eq:reducweakformX0}
a(\bu,\bv)=\ell(\bv)\quad\forall\,\bv\in X_{N_0}.
\end{align}

Set $r(\eta)=p(\eta)/q_A(\eta)$. Let $\mathbb{C}_0^+=\{ r(\eta)\in\mathbb{C}:\Im(r(\eta))\geq 0 \}$, where $\Im(\cdot)$ denotes the imaginary part of a complex number. We now can prove the well-posedness result for problem \eqref{eq:curlcurlequ}.

\begin{lemma}\label{UniqforCurlCurlequ}
For $\bg\in L_t^2(\p\O)$, if $r(\eta)\in\mathbb{C}_0^+$ and $\Re\{q_A(\eta)\}>0$, then \eqref{eq:curlcurlequ} has a unique solution $\bu\in X_{N_0}$ such that
\begin{align}\label{est:u0boundbyf}
\norm{\bu}_{X}\leq C\norm{\bg}_{L^2(\p\O)}
\end{align}
\end{lemma}
\begin{proof}
Define a sesquilinear form $a_+:X_{N_0}\times X_{N_0}\rightarrow\mathbb{C}$ by
\begin{align*}
a_+(\bu,\bv)=q_A(\eta)\left(\bA_0\curl\bu,\curl\bv\right)+\left(\bN_0\bu,\bv\right)-i\left\langle\bu\times\bn,\bv\times\bn\right\rangle.
\end{align*}
Due to $\bA,\,\bN_0\in L^{\infty}(\O)$, it is obvious that $a_+(\cdot,\cdot)$ is continuous on $X_{N_0}\times X_{N_0}$. Taking $\bv=\bu$, we obtain
\begin{align}\label{ineq:apluscoercivity}
\Re\{(1+i) a_+(\bu,\bu) \} &=\Re\{ (1+i)q_A(\eta) \}\left(\bA_0\curl\bu,\curl\bu\right) \notag\\
&\hspace{0.8in}+(\bN_0\bu,\bu)+\norm{\bu\times\bn}^2_{L^2(\p\O)},\notag\\
&\geq \hat{C}\norm{\curl\bu}^2_{L^2(\O)}+(\bN_0\bu,\bu)+\norm{\bu\times\bn}^2_{L^2(\p\O)}, \notag\\
&\geq \alpha\norm{\bu}^2_X,
\end{align}
where we use the positive definiteness of $\bA_0$ and $\bN_0$. Thus, the rotated form $(1+i)a_+(\cdot,\cdot)$ is coercive on $X_{N_0}$. 
By the Lax-Milgram lemma, there exists a unique $\bof\in X_{N_0}$ such that
\begin{align*}
a_+(\bof,\bv)=\ell(\bv) \quad\forall\,\bv\in X_{N_0}.
\end{align*}
Meanwhile, by $|\ell(\bv)|\leq \norm{\bg}_{L^2(\p\O)}\norm{\bv}_X$, we have that
\begin{align}\label{est:fboundbyg}
\norm{\boldsymbol{f}}_X\leq C\norm{\boldsymbol{g}}_{L^2(\partial\Omega)}.
\end{align}
Define a mapping $K:X_{N_0}\rightarrow X_{N_0}$ such that, for each $\bu\in X_{N_0}$, $K\bu\in X_{N_0}$,
\begin{align}
a_+(K\bu,\bv) = -\left(1+p(\eta)\right)(\bN_0\bu,\bv) \quad\forall\,\bv\in X_{N_0}.
\end{align}
Since $H^1_0(\O)\subset S$, every $\bu\in X_{N_0}$ satisfies $\div(\bN_0\bu)=0$ in $\O$. Moreover, $\bu\in H(\text{curl},\O)$ and $\bu\times\bn\in [L^2(\p\O)]^3$. By the weighted Maxwell compactness property, the embedding $X_{N_0}\hookrightarrow [L^2(\O)]^3$ is compact. Hence, the operator $K$ is compact.

It is clear that, if $\bu$ solves~\eqref{eq:reducweakformX0}, $\bu$ satisfies $a_+(\bu+K\bu,\bv)=\ell(\bv)$. Hence, $\bu$ solves the operator equation
\begin{align}\label{eq:FredholmOperFunc}
(I+K)\bu=\bof.
\end{align}
Since $K$ is compact, $I+K$ is a Fredholm operator of index zero. It therefore suffices to prove uniqueness for the corresponding homogeneous problem.
Let $\bg=\boldsymbol{0}$. Taking $\bv=\bu$ gives
\begin{align}\label{eq:homogeneousIKu0}
q_A(\eta)(\bA_0\curl\bu,\curl\bu)-p(\eta)(\bN_0\bu,\bu)-i\norm{\bu\times\bn}^2_{L^2(\p\O)}=0.
\end{align}
Since $\Re\{q_A(\eta)\}>0$, taking the imaginary part yields
\begin{align*}
\Im\left\{r(\eta)\right\}(\bN_0\bu,\bu)+\frac{\Re\{q_A(\eta)\}}{|q_A(\eta)|^2}\norm{\bu\times\bn}^2_{L^2(\p\O)}=0.
\end{align*}
If $\Im(r(\eta))>0$, the uniform positive definiteness of $\bN_0$ immediately implies $\bu=\boldsymbol{0}$. If $\Im\left(r(\eta)\right)=0$, we have $\bu\times\bn=\boldsymbol{0}$ on $\p\O$. The boundary condition in \eqref{def:MaxwellwithAN} then gives $(\bA_0\curl\bu)\times\bn=\boldsymbol{0}$ on $\p\O$. If $r(\eta)<0$, taking the real part of \eqref{eq:homogeneousIKu0} gives $\left( \bA_0\curl\bu,\curl\bu\right)+|r(\eta)|\left(\bN_0\bu,\bu\right)=0,$ and hence $\bu=\boldsymbol{0}$. If $r(\eta)>0$, the conditions $\bu\times\bn=0$ and $(\bA_0\curl\bu)\times\bn=0$ on $\p\O$, as well as the unique continuation property for the Maxwell system~\cite{Monk:2003:FEMMWE}, imply that $\bu=\boldsymbol{0}$ in $\O$. Therefore, the homogeneous problem has only the trivial solution. The Fredholm alternative~\cite{Monk:2003:FEMMWE} and \eqref{est:fboundbyg} imply that \eqref{eq:FredholmOperFunc} has a unique solution $\bu\in X_{N_0}$ satisfying \eqref{est:u0boundbyf}.
%
\end{proof}

Consequently, we define a solution operator
\begin{align*}
S_1(\eta): L_t^2(\p\O)\rightarrow X_{N_0}\quad \text{such\,\,that}\quad \bu:=S_1(\eta)\bg.
\end{align*}

\begin{theorem}\label{thm:holomorphicS1}
The operator $S_1:\mathbb{C}_0^+\rightarrow\mathcal{L}\left(L_t^2(\p\O),X_{N_0}\right)$ is holomorphic.
\end{theorem}
\begin{proof}
Let $\eta\in\mathcal{D}$ and $\delta\eta$ be small enough such that $(\eta+\delta\eta)\in \mathcal{D}$. For a given $\bg\in L_t^2(\p\O)$, let $\bu:=S_1(\eta)\bg$ such that
\begin{align*}
q_A(\eta)\left(\bA_0\curl\bu,\curl\bv\right)-p(\eta)\left(\bN_0\bu,\bv\right)-i\left\langle\bu\times\bn,\bv\times\bn\right\rangle=\left\langle\bg,\bv_T\right\rangle\quad\forall\,\bv\in X_{N_0}
\end{align*}
and $\bw:=S_1(\eta+\delta\eta)\bg$ such that
\begin{align*}
q_A(\eta+\delta\eta)\left(\bA\curl\bw,\curl\bv\right)-p(\eta+\delta\eta)\left(\bN_0\bw,\bv\right)-i\left\langle\bw\times\bn,\bv\times\bn\right\rangle=\left\langle\bg,\bv_T\right\rangle,
\end{align*}
For all $\bv\in X_{X_0}$. Subtracting the above two equations, we obtain
\begin{align}\label{eq:wminusu}
a(\bw-\bu,\bv) = -\delta q_A(\eta)\left(\bA_0\curl\bw,\curl\bv \right)+\delta p(\eta)\left(\bN_0\bw,\bv\right).
\end{align}
Problem \eqref{eq:wminusu} has a unique solution $(\bw-\bu)\in X_{N_0}$. Let $\bpsi_A\,\bpsi_N\in X_{N_0}$ be the solutions of
\begin{align*}
a\left( \bpsi_A,\bv \right) = -\left(\bA_0\curl\bu,\curl\bv\right)\qquad\text{and}\qquad
a\left( \bpsi_N,\bv \right) = \left(\bN_0\bu,\bv\right),
\end{align*}
respectively, for all $\bv\in X_{N_0}$. Let $\bphi:=\delta q_A(\eta)\bpsi_A+\delta p(\eta)\bpsi_N$. Then the error $\be:=\bw-\bu-\bphi$ satisfies
\begin{align*}
a(\be,\bv)= -\delta q_A(\eta)\left( \bA_0\curl(\bw-\bu),\curl\bv \right) + \delta p(\eta)\left(\bN_0(\bw-\bu),\bv\right) \quad\forall\,\bv\in X_{N_0}.
\end{align*}
The right-hand side defines a bounded linear functional on $X_{N_0}$.
Taking $\bv=\be$ in the above equation, by Lemma~\ref{UniqforCurlCurlequ}, we have 
\begin{align}\label{este}
\norm{\be}_{X}\leq C\left( |\delta q_A(\eta)| + |\delta p(\eta)|\right)\norm{\bw-\bu}_{L^2(\O)}.
\end{align}
Similarly, by taking $\bv=\bw-\bu$ in \eqref{eq:wminusu}, and again using Lemma \ref{UniqforCurlCurlequ}, we can obtain that
\begin{align}\label{estwminu}
\norm{\bw-\bu}_{X} &\leq C\left(|\delta q_A(\eta)|+|\delta p(\eta)|\right)\norm{\bw}_{L^2(\O)} \notag\\
&\leq C\left(|\delta q_A(\eta)|+|\delta p(\eta)|\right)|\norm{\bg}_{L^2(\p\O)}.
\end{align}
Combining \eqref{este} and \eqref{estwminu}, we have
\begin{align}
\norm{\be}_{X}\leq C\left(|\delta q_A(\eta)|+|\delta p(\eta)|\right)^2\norm{\bg}_{L^2(\p\O)}.
\end{align}
Assume that $q_A(\eta)$ and $p(\eta)$ are holomorphic at $\eta$. In particular, we have $|\delta q_A(\eta)|+|\delta p(\eta)|\leq C|\delta\eta|$ for sufficiently small $|\delta\eta|$. Hence, we have
\begin{align}
\frac{\norm{\be}_{X}}{|\delta\eta|}\leq C|\delta\eta|\norm{\bg}_{L^2(\p\O)}\rightarrow 0 \quad \text{as} \quad \delta\eta\rightarrow 0.
\end{align}
Therefore, as $\delta\eta\rightarrow 0$, we have
\begin{align*}
&\bigg\|\frac{S_1(\eta+\delta\eta)\bg-S_1(\eta)\bg}{\delta p(\eta)}-q^\prime_A(\eta)\bpsi_A-p^\prime(\eta)\bpsi_N\bigg\|_{X} \\
&\hspace{0.4in}\leq\bigg|\frac{\delta q_A(\eta)}{\delta\eta}-q^\prime_A(\eta)\bigg|\norm{\bpsi_A}_{X}+\bigg|\frac{\delta p(\eta)}{\delta\eta}-p^\prime(\eta) \bigg|\norm{\bpsi_N}_{X}+\frac{\norm{\be}_X}{\delta\eta}\rightarrow 0.
\end{align*}
Hence $S_1(\eta)\bg$ is holomorphic on $\mathbb{C}_0^+$ for every $\bg\in L^2_t(\p\O)$ and thus $S_1(\eta)$ is holomorphic by Theorem 1.7.1 of \cite{gohberg2009holomorphic}.
\end{proof}

Let $X_I:=\{\bv\in X: (\bv,\grad\xi)=0,\,\,\forall\,\xi\in S\}$. Define $S_0: L_t^2(\p\O)\rightarrow X_I$ by letting $\bv:=S_0(\eta)\bg$ be the unique solution of \eqref{eq:curlcurlequ} with $\bA=I$ and $\bN=I$. By Lemma \ref{UniqforCurlCurlequ}, for each $\eta\in \mathbb{C}^+_0$, the mapping $S_j(\eta): L_t^2(\p\O)\rightarrow X$, $j=0$ or $1$, is bounded. By Theorem~\ref{thm:holomorphicS1}, $S_j:\mathbb{C}^+_0\rightarrow\mathcal{L}\left(L_t^2(\p\O),X\right)$, $j=0,1$ are holomorphic. 
We choose an open bounded set $\hat{\mathcal{D}}$ such that $\hat{\mathcal{D}}\Subset\mathbb{C}^+_0$. Then $S_j, j=0,1,$ are holomorphic on $\hat{\mathcal{D}}$.
Define a trace operator
\begin{align*}
\mathcal{I}: X\rightarrow L_t^2(\p\O), \quad \mathcal{I}(\bu):=(\bu\times\bn)|_{\p\O}.
\end{align*}
Then, $\mathcal{I}$ is continuous, i.e., $\norm{\mathcal{I}\bu}_{L^2(\p\O)}\leq  C\norm{\bu}_{X}$. 

Define $T:\hat{\mathcal{D}}\rightarrow \mathcal{L}\left( L_t^2(\p\O), L_t^2(\p\O) \right)$ such that
\begin{align*}
T(\eta):=\mathcal{I}\left[ S_1(\eta)-S_0(\eta) \right].
\end{align*}

\begin{lemma}
The operator function $T(\eta)$ is holomorphic in $\hat{\mathcal{D}}$. 
\end{lemma}
\begin{proof}
Since $S_1(\eta)$ and $S_0(\eta)$ are holomorphic, $T(\eta)$ is holomorphic in $\hat{\mathcal{D}}$. 
\end{proof}


\begin{theorem}
A complex number $\lambda\in\hat{\mathcal{D}}$ is an eigenvalue of $T$ if and only if it is a transmission eigenvalue of \eqref{MTEPA}.
\end{theorem}
\begin{proof}
Let $\lambda$ be an eigenvalue of $T$. Then there exists $\boldsymbol{0}\neq \bg\in L_t^2(\p\O)$ such that
\begin{align*}
T(\lambda)\bg=\mathcal{I}\left[ S_1(\lambda)-S_0(\lambda) \right]\bg=\boldsymbol{0}.
\end{align*}
Let $\bu:=S_1(\lambda)\bg$ and $\bv:=S_0(\lambda)\bg$. Then $\bu$ and $\bv$ solve 
\begin{equation}\label{eq:RP1}
\left\{
\begin{aligned}
\curl\left(\bA\curl\bu\right)-\lambda\bN\bu &=\boldsymbol{0}\quad\text{in}\,\,\O, \\
(\bA\curl\bu)\times\bn-i(\bn\times\bu)\times\bn &=\bg \quad\text{on}\,\,\p\O,
\end{aligned}
\right.
\end{equation}
and 
\begin{equation}\label{eq:RP2}
\left\{
\begin{aligned}
\curl\curl\bv -\lambda\bv &=\boldsymbol{0}\quad\text{in}\,\,\O, \\
(\curl\bv)\times\bn -i(\bn\times\bv)\times\bn &=\bg \quad\text{on}\,\,\p\O,
\end{aligned}
\right.
\end{equation}
respectively. Meanwhile, the boundary conditions yield
\begin{align*}
(\bA\curl\bu)\times\bn - i(\bn\times\bu)\times\bn =(\curl\bv)\times\bn - i(\bn\times\bv)\times\bn \quad\text{on}\,\,\p\O.
\end{align*}
Combining this with $\mathcal{I}\bu=\mathcal{I}\bv$, we obtain that
\begin{align*}
\bu\times\bn=\bv\times\bn \quad \text{on}\,\p\O \quad \text{and} \quad (\bA\curl\bu)\times\bn=(\curl\bv)\times\bn \quad\text{on}\,\p\O.
\end{align*}
Therefore, $(\bu,\bv)$ is a transmission eigenpair associated with $\lambda$. Moreover, $(\bu,\bv)\neq (\boldsymbol{0},\boldsymbol{0})$. Otherwise the boundary conditions in either \eqref{eq:RP1} or \eqref{eq:RP2} would imply $\bg=\boldsymbol{0}$, which contradicts $\bg\neq\boldsymbol{0}$. Thus $\lambda$ is a transmission eigenvalue of \eqref{MTEPA}.

Conversely, suppose that $\lambda$ is a transmission eigenvalue of \eqref{MTEPA} and let $(\bu,\bv)\neq (\boldsymbol{0},\boldsymbol{0})$ be an associated eigenfunction pair satisfying \eqref{MTEPA}. Define 
\begin{align*}
\bg:=(\bA\curl\bu)\times\bn-i(\bn\times\bu)\times\bn=(\curl\bv)\times\bn -i(\bn\times\bv)\times\bn\in L_t^2(\p\O),
\end{align*}
where the equality follows from the conditions in \eqref{MTEPA}. Then $\bu$ and $\bv$ solve the boundary problems \eqref{eq:RP1} and \eqref{eq:RP2} with data $\bg$, respectively. Hence, $\bu=S_1(\lambda)\bg$ and $\bv=S_0(\lambda)(\bg)$. Using the third equation in \eqref{MTEPA}, one has that
\begin{align*}
T(\lambda)\bg=\mathcal{I}\left( S_1(\lambda)\bg-S_0(\lambda)\bg \right)=\mathcal{I}(\bu-\bv)=0.
\end{align*}
If $\bg=\boldsymbol{0}$, both problems \eqref{eq:RP1} and \eqref{eq:RP2} are homogeneous and $\bu=\bv=\boldsymbol{0}$, which contradicts the nontriviality of the transmission eigenvalue pair. Hence, $\bg\neq \boldsymbol{0}$. Thus $T(\lambda)$ has a nontrivial kernel and $\lambda$ is an eigenvalue of $T$. 

The proof is completed.
\end{proof}

We remark that there are other ways to formulate problem \eqref{MTEPA} as an operator eigenvalue problem. For example, define the space
\begin{align*}
\hat{X}=TH^{-1/2}(\text{div}_{\p\O},\p\O):=\{ \bv\in [H^{-1/2}(\p\O)]^3: \bv\cdot\bn=0,\,\nabla_{\p\O}\cdot\bv\in H^{-1/2}(\p\O) \}. 
\end{align*}
For $\eta\in\mathbb{C}$ not belonging to the interior Dirichlet spectra of 
\begin{align*}
\curl(\bA\curl\bu)-\eta\bN\bu=\boldsymbol{0}\quad\text{and}\quad \curl\curl\bv-\eta\bv=\boldsymbol{0} \quad \text{in } \O,
\end{align*}
define $\hat{\bu}$ and $\hat{\bv}$ as the respective solutions to 
\begin{align*}
\left\{
\begin{aligned}
\curl(\bA\curl\hat{\bu})-\eta\bN\hat{\bu} &=\boldsymbol{0}\quad\text{in}\,\,\O \\
\hat{\bu}\times\bn &=\bg\quad\text{on}\,\,\p\O,
\end{aligned}
\right. 
\quad\text{and}\quad 
\left\{
\begin{aligned}
\curl(\curl\hat{\bv})-\eta\hat{\bv} &=\boldsymbol{0}\quad\text{in}\,\,\O \\
\hat{\bv}\times\bn &=\bg\quad\text{on}\,\,\p\O.
\end{aligned}
\right.
\end{align*}
The Dirichlet-to-Neumann operators $\hat{S}_1:L^2_t(\p\O)\to \hat{X}$ and $\hat{S}_0:L^2_t(\p\O)\to \hat{X}$ associated with the above two systems are respectively defined by 
\begin{align*}
\hat{S}_1(\eta)\bg=(\bA\curl\hat{\bu})\times\bn \quad\text{and}\quad \hat{S}_0(\eta)\bg=(\curl\hat{\bv})\times\bn.
\end{align*}
Then the operator eigenvalue problem is to find $\eta$ and $\bg\neq\boldsymbol{0}$ such that
\begin{align*}
\hat{T}(\eta)\bg:=\left( \hat{S}_1(\eta)-\hat{S}_0(\eta)\right)\bg=\boldsymbol{0}.
\end{align*}



\section{Finite Element Approximation}\label{sec:FEP} 
%
In this section, we propose a finite element approximation for $T(\eta)$. Let $\cT_h$ be a regular tetrahedral mesh of $\O$ with mesh size $h$. Let $V_h\subset H(\text{curl},\O)$ be the N\'ed\'elec's edge finite element space associated with $\cT_h$. Define $V_h^B:=\{ \gamma_T\bv_h|_{\p\O}:\,\bv_h\in V_h \}\subset L^2_t(\p\O)$ as the restriction of $V_h$ on $\p\O$. 
Let $p_h: L_t^2(\p\O)\rightarrow V_h^B$ be the $L^2$-projection defined by
\begin{align}\label{DefPh}
\left\langle\bg-p_h\bg,\bmu_h\right\rangle = 0, \qquad\forall\,\bmu_h\in V_h^B.
\end{align} 
For $\eta\in \hat{\mathbb{C}}$, the finite element formulation for \eqref{eq:curlcurlequ} is to find $\bu_{1,h}\in V_h$ such that
\begin{align}\label{FEMScheme}
a\left(\bu_{1,h},\bv_h\right)=\left\langle p_h\bg, \gamma_T\bv_h\right\rangle\quad\forall\,\bv_h\in V_h,
\end{align} 
where $a(\cdot,\cdot)$ is given in \eqref{eq:defa}. 
%
Since $\gamma_T\bv_h\in V_h^B$, the Galerkin orthogonality preserves, i.e., $a(\bu-\bu_{1,h},\bv_h)=0$ for any $\bv_h\in V_h$.

\noindent\textbf{Assumption A}. For every $r(\eta)\in\mathbb{C}^+_0$, the solution $\bu=S_1(\eta)\bg$ to \eqref{eq:curlcurlequ} satisfies $\bu\in [H^s(\O)]^3$ and $\curl\bu\in [H^s(\O)]^3$ with $1+\delta\leq s\leq r$, $\delta>0$. Furthermore,
\begin{align}\label{assump:regforu}
\norm{\bu}_{H^s(\O)}+\norm{\curl\bu}_{H^s(\O)}\leq C\norm{\bg}_{L^2(\p\O)}.
\end{align}
We also assume the corresponding adjoint regularity uniformly for $r(\eta)\in\mathbb{C}^+_0$, i.e., if $\bw\in X$ solves $a(\bv,\bw)=(\bv,\bof),\,\,\forall\,\bv\in X$, then
\begin{align}\label{assump:regforw}
\norm{\bw}_{H^s(\O)}+\norm{\curl\bw}_{H^s(\O)}\leq C\norm{\bof}_{L^2(\p\O)}.
\end{align}
The same regularity assumptions are made for the problem associated with $S_0(\eta)$. 

\begin{lemma}
Let $r(\eta)\in\mathbb{C}^+_0$ and $\Re\{q_A(\eta)\}>0$. There exists a unique solution $\bu_{1,h}$ to \eqref{FEMScheme}. Moreover,
\begin{align}\label{disStabilityEst}
\norm{\bu_{1,h}}_X\leq C\norm{\bg}_{L^2(\p\O)}.
\end{align}
\end{lemma}
\begin{proof}
Using Cauchy-Schwarz inequality, it is easy to see that the sesquilinear form $a(\cdot,\cdot)$ is bounded on $X\times X$, i.e., there exists $C_1>0$ such that
\begin{align}\label{ineq:boundofa}
|a(\bu,\bv)|\leq C_1\norm{\bu}_X\norm{\bv}_X \quad\forall\,\bu,\,\bv\in X,
\end{align}
where we can take $C_1=\max\{|q_A(\eta)|\norm{\bA_0}_{L^\infty(\O)}, |p(\eta)|\,\norm{\bN_0}_{L^\infty(\O)}, 1\}$. Moreover, we choose $\sigma>\Re\{p(\eta)\}-\Im \{p(\eta)\}$, then for any $\bv\in X$, we have
\begin{align*}
&\Re\left\{(1+i)\left[ a(\bv,\bv) +\sigma(\bN\bv,\bv)\right] \right\} \notag\\
&\hspace{0.4in}=\Re\{(1+i)q_A(\eta)\}(\bA_0\curl\bv,\curl\bv) \\
&\hspace{0.8in}+\left(\sigma-\Re\{p(\eta)\}+\Im\{p(\eta)\}\right)(\bN\bv,\bv)+\norm{\bv\times\bn}^2_{L^2(\p\O)} \notag\\
&\hspace{0.4in}\geq C_A\norm{\curl\bv}^2_{L^2(\O)}+\left(\sigma-\Re\{p(\eta)\}+\Im \{p(\eta)\}\right)N_*\norm{\bv}^2_{L^2(\O)} + \norm{\bv\times\bn}^2_{L^2(\p\O)} \notag\\
&\hspace{0.4in}\geq C_2\norm{\bv}^2_X,
\end{align*}
where $C_2=\min\{ A_*, (\sigma-\Re p(\eta)+\Im p(\eta))N_*, 1 \}$. Therefore, for any $\bv\in X$,
\begin{align}\label{ineq:gradingofa}
\norm{\bv}_X^2\leq C_0\big(|a(\bv,\bv)|+\norm{\bv}^2_{L^2(\O)}\big),
\end{align}
where $C_0>0$ are independent of $h$. Using \eqref{ineq:gradingofa}, the continuous uniqueness, the standard Schatz argument implies that the Galerkin problem is uniquely solvable for sufficiently small $h$, with the uniform stability estimate \eqref{disStabilityEst}.
\end{proof}

For $r(\eta)\in\mathbb{C}^+_0$ and $\bg\in L_t^2(\p\O)$, define the discrete solution operator $S_{1,h}(\eta): L_t^2(\p\O)\rightarrow V_h$ by 
\begin{align*}
\bu_{1,h}=S_{1,h}(\eta)\bg,
\end{align*}
where $\bu_{1,h}$ is the unique solution of \eqref{FEMScheme}. Let $\mathcal{I}_h: V_h\rightarrow V_h^B$ be the discrete tangential trace map defined by $\mathcal{I}_h\bu_{1,h}:=(\bu_{1,h}\times\bn)|_{\p\O}$. Then the discrete boundary trace corresponding to $\bg$ is defined by
\begin{align}
\bof_{1,h}:=\mathcal{I}_hS_{1,h}(\eta)\bg\in V_h^B, \qquad \text{i.e.,}\,\,\bof_{1,h}=(\bu_{1,h}\times\bn)|_{\p\O}.
\end{align}

Let $\Pi_h$ denote the standard $H(\text{curl})$-conforming interpolation. The following approximation properties \cite{Nedelec:1980:EdgeElements, Nedelec:1986:EdgeElements} hold.
\begin{lemma}\label{InterEst}
If $\bv\in [H^s(\O)]^3$ and $\curl\bv\in [H^s(\O)]^3$ for $\frac{1}{2}+\delta\leq s\leq r$ and $\delta>0$, on each $T\in\cT_h$ it holds that
\begin{align*}
\norm{\bv-\Pi_h\bv}_{L^2(T)}+\norm{\curl(\bv-\Pi_h\bv)}_{L^2(T)} &\leq Ch^s_T\left(\norm{\bv}_{H^s(T)}+\norm{\curl\bv}_{H^s(T)} \right) \\
\norm{\curl(\bv-\Pi_h\bv)}_{L^2(T)} &\leq Ch^s_T\norm{\curl\bv}_{H^s(T)},
\end{align*}
where $C>0$ is a constant that is independent of the mesh size.
\end{lemma}



We now provide the convergent result for the solution on boundary as follows.
\begin{theorem}\label{fferrest}
Let $\eta\in\mathbb{C}_0^+$ be fixed and $\bu$ be the solution of \eqref{eq:curlcurlequ} satisfying the Assumption A. Let $\bof_1=\mathcal{I}S_1(\eta)\bg$ and $\bof_{1,h}=\mathcal{I}_hS_{1,h}(\eta)\bg$ with $\bg\in L_t^2(\p\O)$. Then
\begin{align*}
\norm{\bof_{1,h}-\bof_1}_{L^2(\p\O)} \leq Ch^{s-\frac{1}{2}}\norm{\bg}_{L^2(\p\O)},
\end{align*}
where $C>0$ is independent of the mesh size $h$ and $\bg$. 
\end{theorem}
\begin{proof}
Let $\eta\in\mathbb C_0^+$ be fixed. Let $\bu=S_1(\eta)\bg\in X$ and $\bu_{1,h}=S_{1,h}(\eta)\bg\in V_h$ be the solutions of \eqref{eq:curlcurlequ} and \eqref{FEMScheme}, respectively. By the definition of $p_h$ in \eqref{DefPh}, the discrete right-hand side preserves Galerkin orthogonality, namely,
\begin{align}\label{eq:GlkOrth}
a(\bu-\bu_{1,h},\bv_h)=0 \quad\forall\,\bv_h\in V_h.
\end{align}
For any $\bv_h\in V_h$, by using \eqref{eq:GlkOrth}, 
\eqref{ineq:boundofa} and \eqref{ineq:gradingofa}, we obtain
\begin{align}\label{ineq:Xpolluted}
\norm{\bu-\bu_{1,h}}_X^2 &\leq C\left( |a(\bu-\bu_{1,h},\bu-\bu_{1,h})| +\norm{\bu-\bu_{1,h}}_{L^2(\O)}^2\right)  \\
&=C\left( |a(\bu-\bu_{1,h},\bu-\bv_h)| +\norm{\bu-\bu_{1,h}}_{L^2(\O)}^2\right) \nonumber \\
&\leq C\left(\norm{\bu-\bu_{1,h}}_X\norm{\bu-\bv_h}_X +\norm{\bu-\bu_{1,h}}_{L^2(\O)}^2\right). \notag
\end{align}
Let $\bw\in X$ solve the adjoint problem
\begin{align}\label{eq:dual}
a(\bv,\bw)=(\bv,\bu-\bu_{1,h}) \qquad\forall\,\bv\in X.
\end{align}
By \eqref{assump:regforw} in Assumption A, we obtain
\begin{equation}\label{eq:adj_reg}
\norm{\bw}_{H^s(\O)}+\norm{\curl \bw}_{H^s(\O)}\le C\norm{\bu-\bu_{1,h}}_{L^2(\O)}.
\end{equation}
For any $\bw_h\in V_h$, taking $\bv=\bu-\bu_{1,h}$ in \eqref{eq:dual}, then \eqref{eq:GlkOrth} and \eqref{ineq:boundofa} gives
\begin{align*}
\norm{\bu-\bu_{1,h}}_{L^2(\O)}^2 &=a(\bu-\bu_{1,h},\bw)=a(\bu-\bu_{1,h},\bw-\bw_h) \\
&\leq C_1\norm{\bu-\bu_{1,h}}_X\norm{\bw-\bw_h}_X\leq Ch^s\norm{\bu-\bu_{1,h}}_X\norm{\bu-\bu_{1,h}}_{L^2(\O)},
\end{align*}
where we have used the regularity of the solution for the adjoint problem~\eqref{eq:dual}. Hence, we obtain
\begin{align}\label{eq:L2_vs_X}
\norm{\bu-\bu_{1,h}}_{L^2(\O)}\leq Ch^s\norm{\bu-\bu_{1,h}}_X.
\end{align}
Substituting \eqref{eq:L2_vs_X} into \eqref{ineq:Xpolluted} yields
\begin{align*}
\norm{\bu-\bu_{1,h}}_X\leq C\norm{\bu-\bv_h}_X+Ch^{2s}\norm{\bu-\bu_{1,h}}_X.
\end{align*}
For $h$ sufficiently small, the last term can be absorbed, and hence 
\begin{align}
\norm{\bu-\bu_{1,h}}_X\leq C\inf_{\bv\in V_h}\norm{\bu-\bv_h}_X \quad\forall\,\bv\in V_h.
\end{align}
Choosing $\bv_h=\Pi_h\bu$ and applying Lemma~\ref{InterEst} together with regularity assumption \eqref{assump:regforu}, we conclude
\begin{align}\label{eq:X_rate}
\norm{\bu-\bu_{1,h}}_X\leq Ch^s\left(\norm{\bu}_{H^s(\O)}+\norm{\curl\bu}_{H^s(\O)}\right)\leq Ch^s\norm{\bg}_{L^2(\p\O)},
\end{align}
which combines with \eqref{eq:L2_vs_X} yields
\begin{align}
\norm{\bu-\bu_{1,h}}_{L^2(\O)}\leq Ch^{2s}\norm{\bg}_{L^2(\p\O)}.
\end{align}
For the boundary error, by triangle inequality, we have
\begin{align}\label{ineq:bd_split}
\norm{(\bu-\bu_{1,h})\times\bn}_{L^2(\p\O)}\leq \norm{(\bu-\Pi_h\bu)\times\bn}_{L^2(\p\O)}+\norm{(\Pi_h\bu-\bu_{1,h})\times\bn}_{L^2(\p\O)}.
\end{align}
For the term $\norm{(\bu-\Pi_h\bu)\times\bn}_{L^2(\p\O)}$, using the trace theorem (cf. Lemma 5.52 of \cite{Monk:2003:FEMMWE}) and Lemma \ref{InterEst}, we have
\begin{align*}
\norm{(\bu-\Pi_h\bu)\times\bn}_{L^2(\p\O)}\leq Ch^{s-\frac{1}{2}}\left(\norm{\bu}_{H^s(\O)}+\norm{\curl\bu}_{H^s(\O)}\right)\leq Ch^{s-\frac{1}{2}}\norm{\bg}_{L^2(\p\O)}.
\end{align*}
For the term $\norm{(\Pi_h\bu-\bu_{1,h})\times\bn}_{L^2(\p\O)}$, we employ the standard discrete trace inequality for N\'ed\'elec edge element, i.e.,
\begin{align*}
&\norm{(\Pi_h\bu-\bu_{1,h})\times\bn}_{L^2(\p\O)} \\
&\hspace{0.4in}\leq C\left(h^{-1/2}\norm{\Pi_h\bu-\bu_{1,h}}_{L^2(\O)}+h^{1/2}\norm{\curl(\Pi_h\bu-\bu_{1,h})}_{L^2(\O)}\right) \\
&\hspace{0.4in}\leq Ch^{-1/2}\norm{\Pi_h\bu-\bu_{1,h}}_X \\
&\hspace{0.4in}\leq Ch^{-1/2}\left(\norm{\Pi_h\bu-\bu}_X+\norm{\bu-\bu_{1,h}}_X\right) \leq  Ch^{s-\frac{1}{2}}\norm{\bg}_{L^2(\p\O)}.
\end{align*}
Combining the above two inequalities with \eqref{ineq:bd_split} yields
\begin{align*}
\norm{\bof_{1,h}-\bof_{1}}_{L^2(\p\O)}\leq Ch^{s-\frac{1}{2}}\norm{\bg}_{L^2(\p\O)}.
\end{align*}
This completes the proof.

\end{proof}

Set $\bA=I$ and $\bN=I$. Consider the problem of finding $\bu_{0,h}\in V_h$ such that 
\begin{align}\label{FEMScheme2}
a(\bu_{0,h},\bv_h) = \left\langle p_h\bg, \gamma_T\bv_h\right\rangle \qquad\forall\,\bv_h\in V_h.
\end{align}
Similarly, we define the solution operator $S_{0,h}(\eta): L_t^2(\p\O)\rightarrow V_h$ for \eqref{FEMScheme2} associated with \eqref{FEMScheme2} by $\bu_{0,h}=S_{0,h}(\eta)\bg$ and 
\begin{align}
\bof_{0,h}=\mathcal{I}S_{0,h}(\eta)\bg:=(\bu_{0,h}\times\bn)|_{\p\O}.
\end{align}
From Theorem \ref{fferrest}, one has that 
\begin{align}
\norm{\bof_0-\bof_{0,h}}_{L^2(\p\O)}\le Ch^{s-\frac{1}{2}}\norm{\bg}_{L^2(\p\O)} \qquad\forall\,\bg\in L_t^2(\p\O),
\end{align}
where $\bof_0=\mathcal{I}S_0(\eta)\bg$.
Let $T_h(\eta)$ be a finite element approximation of $T(\eta)$, defined by 
\begin{align}\label{eq:defTh}
T_h(\eta):=\mathcal{I}_h\left( S_{1,h}(\eta)-S_{0,h}(\eta) \right),
\end{align}
where $\mathcal{I}_h: V_h\rightarrow V_h^B$ is the restriction operator. 

\begin{lemma}\label{lemma:convph}
For every $\bg\in L^2_t(\p\O)$, the projection $p_h\in \mathcal{L}\left( L_t^2(\p\O), V_h^B \right)$ defined in \eqref{DefPh} satisfies
\begin{align}\label{conphgtog}
\lim_{h\rightarrow 0}\norm{p_h\bg}_{L^2(\p\O)}=\norm{\bg}_{L^2(\p\O)}.
\end{align}
\end{lemma}
\begin{proof}
By the definition of \eqref{DefPh}, we have that 
\begin{align}\label{eq:bestapprox_ph}
\norm{\bg-p_h\bg}_{L^2(\p\O)} =\inf_{\bv_h\in V_h^B}\norm{\bg-\bv_h}_{L^2(\p\O)}.
\end{align}
We first let $\bg=(\bu\times\bn)|_{\p\O}$ with $\bu\in [H^s(\O)]^3$, $s>1$. Let $\bv_h=(\Pi_h\bu\times\bn)|_{\p\O}$, by using the trace inequality (cf. Section 5.6 in \cite{Monk:2003:FEMMWE}) and Lemma~\ref{InterEst}, we obtain
\begin{align}\label{ineq:convgminphg}
\norm{\bg-p_h\bg}_{L^2_t(\p\O)} &\leq\norm{(\bu-\Pi_h\bu)\times\bn}_{L^2(\p\O)}\\
&\leq C\sum_{T\in\cT_h^\p}\left( h_T^{-1/2}\norm{\bu-\Pi_h\bu}_{L^2(T)}+h_T^{1/2}|\bu-\Pi_h\bu|_{H^1(T)}\right), \notag\\
&\leq C\sum_{T\in\cT_h^\p}h_T^{s-\frac{1}{2}}\norm{\bu}_{H^s(T)}\leq Ch^{s-\frac{1}{2}}\norm{\bu}_{H^s(\O)}, \notag
\end{align}
where $\cT_h^\p:=\{ T\in\cT_h:\p T\cap\p\O\neq\emptyset\}$. Hence, $\norm{\bg-p_h\bg}_{L^2(\p\O)}\rightarrow 0$ as $h\rightarrow 0$. We now consider an arbitrary $\bg\in L^2_t(\p\O)$. By the density of smooth tangential fields in $L^2_t(\p\O)$, for any $\varepsilon>0$, there exists a smooth tangential field $\bg_\varepsilon$ of the above form such that $\norm{\bg-\bg_\varepsilon}_{L^2(\p\O)}<\varepsilon$. Therefore,
\begin{align*}
\norm{\bg-p_h\bg}_{L^2(\p\O)} &\leq \norm{\bg-\bg_\varepsilon}_{L^2(\p\O)} +\norm{\bg_\varepsilon-p_h\bg}_{L^2(\p\O)} + \norm{p_h(\bg_\varepsilon-\bg)}_{L^2(\p\O)} \\
&\leq 2\varepsilon +\norm{\bg_\varepsilon-p_h\bg_\varepsilon}_{L^2(\p\O)}. 
\end{align*}
For fixed $\bg_\varepsilon$, the last term tends to zero as $h\rightarrow 0$. Thus 
\begin{align*}
\lim_{h\rightarrow 0}\sup\norm{\bg-p_h\bg}_{L^2(\p\O)}\leq 2\varepsilon.
\end{align*} 
Since $\varepsilon>0$ is arbitrary, $\norm{\bg-p_h\bg}_{L^2(\p\O)}\rightarrow 0$ as $h\rightarrow 0$. 
Therefore, for every $\bg\in L^2_t(\p\O)$, we have
\begin{align*}
\norm{\bg}^2_{L^2(\p\O)}=\norm{\bg-p_h\bg}^2_{L^2(\p\O)}+\norm{p_h\bg}^2_{L^2(\p\O)},
\end{align*}
which leads to \eqref{conphgtog}.
\end{proof}

\begin{lemma}\label{lemma:conTh}
There exists $h_0>0$ such that, for every compact set $\mathcal{D}\Subset \mathbb{C}_0^+$, 
\begin{align}
\sup_{0<h<h_0}\sup_{\eta\in\mathcal D}\norm{T_h(\eta)}_{\mathcal L(L^2_t(\p\O),\,L_t^2(\p\O))}<\infty.
\end{align}
\end{lemma}
\begin{proof}
Fix $\eta\in\mathcal D,\,(0<h<h_0)$ and $\bg\in L_t^2(\p\O)$. Let $\bu_{1,h}=S_{1,h}(\eta)\bg$, $\bu_{0,h}=S_{0,h}(\eta)\bg$. By \eqref{disStabilityEst}, we have $\norm{S_{1,h}(\eta)\bg}_X\leq C\norm{\bg}_{L^2(\p\O)}$. 
Since
\begin{align*}
\norm{\mathcal I_h\bv_h}_{L^2(\p\O)}=\norm{\bv_h\times\bn}_{L^2(\p\O)}\leq\norm{\bv_h}_X,
\end{align*}
we obtain $\norm{\mathcal I_hS_{1,h}(\eta)\bg}_{L^2(\p\O)}\leq C\norm{\bg}_{L^2(\p\O)}.$
Similarly, we get $\norm{\mathcal I_hS_{0,h}(\eta)\bg}_{L^2(\p\O)}\leq C\norm{\bg}_{L^2(\p\O)}$. 
Consequently, 
\begin{align*}
\norm{T_h(\eta)}=\norm{\mathcal{I}_h\left( S_{1,h}(\eta)-S_{0,h}(\eta) \right)}\leq \norm{\mathcal{I}_hS_{1,h}(\eta)}+\norm{\mathcal{I}_hS_{0,h}(\eta)} \leq \tilde{C}\norm{\bg}_{L^2(\p\O)}.
\end{align*}
The bound is uniform with respect to $0<h<h_0$ and $\eta\in \mathcal D$. The proof is complete.
\end{proof}

\begin{lemma}\label{lemma:conThph}
Let $\eta\in\mathcal D\Subset\mathbb{C}_0^+$, then for all $\bg\in L^2_t(\p\O)$ we have
\begin{align}\label{conThphg}
\lim_{h\rightarrow 0}\norm{T_h(\eta)p_h\bg-p_hT(\eta)\bg}_{L^2(\p\O)}=0.
\end{align}
\end{lemma}
\begin{proof}
Using the triangle inequality, Theorem \ref{fferrest}, Lemma \ref{lemma:convph} we have
\begin{align*}
&\norm{T_h(\eta)p_h\bg-p_hT(\eta)\bg}_{L^2(\p\O)} \\
&\hspace{0.6in}\leq \norm{T_h(\eta)p_h\bg-T(\eta)\bg}_{L^2(\partial\Omega)}+\norm{(I-p_h)T(\eta)\bg}_{L^2(\partial\Omega)}, \\
&\hspace{0.6in}\leq Ch^{s-\frac{1}{2}}\norm{\bg}_{L^2(\p\O)}.
\end{align*}
This completes the proof. 
\end{proof}

Now we are ready to present the main convergence result.
\begin{theorem}
Let $\lambda\in \sigma(T)$ and $h$ be small enough. Assume that $T(\eta)\in\mathcal{L}(\hat{X},\hat{Y})$ is Fredholm of index 0 for $r(\eta)\in\mathbb{C}_0^+$, where $\hat{X}$ and $\hat{Y}$ are Sobolev spaces satisfying $L_t^2(\p\O)\subset \hat{X}$ and $\hat{Y}\subset L_t^2(\p\O)$.  There exist $\lambda_h\in\sigma(T_h)$ such that $\lambda_h\rightarrow \lambda$ as $h\rightarrow 0$. For any sequence $\lambda_h\in \sigma(T_h)$ converging to $\lambda$ and any associated eigenvectors
\begin{align*}
\bg_{0,h}\in\mathcal{N}\left( T_h(\lambda_h) \right),\qquad \norm{\bg_{0,h}}_{L^2(\p\O)}=1,
\end{align*}
the following estimates hold:
\begin{align*}
|\lambda_h-\lambda| &\leq Ch^{(s-\frac{1}{2})/r_0} \\
\inf_{\bg\in\mathcal{N}(T(\lambda))}\norm{\bg_{0,h}-p_h\bg_0}_{L^2(\p\O)} &\leq Ch^{(s-\frac{1}{2})/r_0},
\end{align*}
where $r_0:=\kappa(T,\lambda)$.
\end{theorem}
\begin{proof}

Fix an eigenvalue $\lambda\in\sigma(T)$. Let $\mathcal D\subset\mathbb{C}_0^+$ be a compact set such that $\lambda\in\mathcal{D}$ and $\mathcal{D}$ contains no other spectral points of $T$ in a small neighborhood of $\lambda$.
Let $\{ h_n\}$ be any sequence with $h_h\to0$ as $n\rightarrow \infty$, and set $T_n(\eta):=T_{h_n}(\eta)$, $X=Y=L^2_t(\p\O)$, $X_n=Y_n:=V_{h_n}^B$, $p_n=q_n:=p_{h_n}$. We verify assumptions (a1)--(a4) of Theorem \ref{thm:abstractapprholoFredmOpt}.

%

Assumption (a1)follows from the construction $T_n(\eta)$. For each fixed $n$, $T_n(\eta)$ acts on the finite-dimensional space $V_{n}^B$, hence it is Fredholm of index zero and holomorphic on $\mathcal D$. Assumption (a2) is exactly Lemma \ref{lemma:conTh}. Due to Lemma \ref{lemma:conThph}, assumption (a3) holds with $q_n=p_n$. Moreover, by Theorem \ref{fferrest}, there exist $\delta>0$ and $C>0$ such that
\begin{equation}\label{eq:eps_rate}
\varepsilon_h:=\max_{|\eta-\lambda|\le\delta}\ \max_{\bg\in G(\lambda),\,\|\bg\|_X=1}
\|T_h(\eta)p_h \bg-p_hT(\eta)\bg\|_{L^2(\p\O)}\le C\,h^{\,s-\frac12}.
\end{equation}


We now verify assumption (a4). Let $N^\prime\subset\mathbb{N}$ and $\bv_n\in V_n^B$, $n\in N^\prime$ be a subsequence with $\norm{\bv_n}_{L^2(\p\O)}$ bounded and 
\begin{align}\label{eq:A3b_ass}
\lim_{n\rightarrow\infty}\norm{T_n(\lambda)\bv_h-p_n\boldsymbol{y}}_{L^2(\p\O)} \qquad\text{for some } \boldsymbol{y}\in L_t^2(\p\O).
\end{align}
We now show that there exists a subsequence $N^{\prime\prime}\subset N^\prime$ and a $\bv\in L_t^2(\p\O)$ such that $\lim_{N^{\prime\prime}\ni n\to\infty}\norm{\bv_n-p_n\bv}_{L^2(\p\O)}=0$.

If $\lambda\in\rho(T)$, then $T(\lambda)^{-1}$ exists and is bounded on $L^2(\p\O)$. Set $\bv:=T(\lambda)^{-1}\boldsymbol{y}\in L^2(\p\O)$. By Lemma \ref{lemma:conThph} with $\eta=\lambda$, we have
\begin{align*}
\norm{T_n(\lambda)p_n\bv-p_nT(\lambda)\bv}_{L^2(\p\O)}=\norm{T_n(\lambda)p_n\bv-p_n\boldsymbol{y}}_{L^2(\p\O)}\to0.
\end{align*}
Combining with \eqref{eq:A3b_ass} yields $\norm{T_n(\lambda)(\bv_n-p_n\bv)}_{L^2(\p\O)}\to0$. Since $\lambda\in\rho(T)$ and $\{ T_n(\lambda) \}$ converges regularly to $T(\lambda)$, there exists $n_0$ such that for all $n\geq n_0$, $T_n(\lambda)$ is invertible and $\sup_{n\geq n_0}\norm{T_n(\lambda)^{-1}}<\infty$. Therefore,
\begin{align*}
\norm{\bv_n-p_n\bv}_{L^2(\p\O)}\leq \norm{T_n(\lambda)^{-1}}\,\norm{T_n(\lambda)(\bv_n-p_n\bv)}_{L^2(\p\O)}\to0,
\end{align*}
which proves the assumption (a4).

If $\lambda\in\sigma(T)$, let $G(\lambda)$ be the generalized eigenspace associated with $\lambda$. Then $T(\lambda)$ induces a bounded invertible map from $R(T(\lambda))$ onto the quotient space $L^2_t(\p\O)/G(\lambda)$. Assume $\boldsymbol{y}\in R(T(\lambda))$ and define $\bv^\prime\in L^2(\p\O)/G(\lambda)$ by $\bv^\prime:=T(\lambda)^{-1}\boldsymbol{y}$. Choose a representative in $L^2_t(\p\O)$. Define the cosets $\hat{\bv}_n:=(\bv_n-p_n\bv^\prime)|_{G(\lambda)}\in L^2(\p\O)/G(\lambda)$. Since $G(\lambda)$ is finite dimensional, the unit ball of the quotient space is sequentially compact up to the extracting a subsequence. Hence, $\{\hat{\bv}_n \}$ admits a convergent subsequence with $\hat{\bv}_n\to\hat{\bv}$ in $L^2_t(\p\O)/G(\lambda)$. Let $\bv\in L^2_t(\p\O)$ be any representative of the limit coset $\hat{\bv}+\bv^\prime|_{G(\lambda)}$. Then, by construction, $\norm{(\bv_n-p_n\bv)|_{G(\lambda)}}_{L^2(\p\O)/G(\lambda)}\to0$. 
In both cases, assumption (a4) holds.



All assumptions (a1)--(a4) of Theorem \ref{thm:abstractapprholoFredmOpt} are verified. Hence, for $\lambda\in\sigma(T)$, there exists $h_0>0$ such that for all $h<h_0$, there is $\lambda_h\in\sigma(T_n)$ with $\lambda_h\to\lambda$ as $h\to0$. Moreover, Theorem \ref{thm:abstractapprholoFredmOpt} yields
\begin{equation}\label{eq:abs_est}
|\lambda_h-\lambda|\le C\,\varepsilon_h^{\,1/\kappa(T,\lambda)},\qquad
\inf_{\bg\in\mathcal N(T(\lambda))}\|\bg_{0,h}-p_h \bg\|_{L^2(\p\O)}\le C\,\varepsilon_h^{\,1/\kappa(T,\lambda)},
\end{equation}
Finally, using \eqref{eq:eps_rate}, we obtain
\[
|\lambda_h-\lambda|\le C\,h^{(s-\frac12)/\kappa(T,\lambda)},\qquad
\inf_{\bg\in\mathcal N(T(\lambda))}\|\bg_{0,h}-p_h \bg\|_{X_h}\le C\,h^{(s-\frac12)/\kappa(T,\lambda)}.
\]
This completes the proof.
\end{proof}

For simple eigenvalues, the following corollary holds.
\begin{corollary}
For a simple eigenvalue $\lambda$, there exist $\lambda_h\in\sigma(T_h)$ such that
\begin{align*}
|\lambda_h-\lambda|\leq Ch^{s-\frac{1}{2}}.
\end{align*}
\end{corollary}

\section{Numerical Experiments}\label{sec:pmCIMandNumExp}    

In this section, we first provide the computational procedure for the proposed method in \cref{sec:FEP} and then present numerical results to demonstrate the effectiveness of this method. 
%

\subsection{A Parallel Multi-step Contour Integral Method}\label{sec:pmCIM} 

We use a parallel multi-step contour integral method (pmCIM)~\cite{XiSun2024pmCIM} for computing eigenvalues of $T_h$. This method consists of three stages, i.e., spectral-indicator screening, local eigenvalue extraction by Beyn's method, and a final verification. The theoretical background of the spectral indicator method (SIM) and Beyn's method has been extensively discussed in~\cite{beyn2012integral, XiSun2024pmCIM}. 

Let $\mathcal D\subset \mathbb{C}$ be a bounded simply connected region. Denote by $\Gamma$ the positively oriented circle that circumscribes $\mathcal D$. Assume that $T_h(\eta)$ is invertible for all $\eta\in\Gamma$. The spectral projection operator is defined by
\begin{align}\label{eq:proj}
P_{\Gamma}:=\frac{1}{2\pi i}\int_{\Gamma}T_h(\eta)^{-1}\,\mathrm d\eta.
\end{align}
If $\Gamma$ encloses no eigenvalues of $T_h$, then $P_{\Gamma}\bv_h=0$ for any $\bv_h\in V_h^B$. Otherwise, $P_\Gamma\bv_h\neq 0$. 
%
%
The construction of $T_h(\eta)$ is based on the following finite element matrices. Let $\{\boldsymbol\phi_i\}_{i=1}^{N}$ be the N\'ed\'elec basis of $V_h$ and let $\{\boldsymbol\psi_j\}_{j=1}^{N_B}$ be the induced basis of $V_h^B$ on $\partial\Omega$. The matrices used in the discretization are summarized in Table~\ref{tab:matrix_assembly}.

\begin{table}[htbp]
\small
\caption{Finite element matrices used in the construction of $T_h(\eta)$.}\label{tab:matrix_assembly}
\begin{center}
  \begin{tabular}{c c c p{5.0cm}} \hline
   \toprule
Matrix & Size & Definition & Construction \\
\midrule

$A_h^1$ & $N\times N$ & $\left(\bA\curl\bu_h,\curl\bv_h\right)$ & $(A_h^1)_{ij} = (\boldsymbol A\nabla\times\boldsymbol\phi_j,\nabla\times\boldsymbol\phi_i)_{\Omega}$ \\

$M_h^B$ & $N\times N$ & $\left\langle \bu_h\times\bn,\bv_h\times\bn \right\rangle$ & $(M_h^B)_{ij} = \langle(\boldsymbol\phi_j)_T, (\boldsymbol\phi_i)_T \rangle_{\partial\Omega}$ \\

$M_h^N$ & $N\times N$ & $\left(\bN\bu_h,\bv_h\right)$ & $(M_h^N)_{ij} = (\boldsymbol N\boldsymbol\phi_j, \boldsymbol\phi_i)_{\Omega}$ \\

$A_h^0$ & $N\times N$ & $\left(\curl\bu_h,\curl\bv_h\right)$ & $(A_h^0)_{ij} = (\nabla\times\boldsymbol\phi_j, \nabla\times\boldsymbol\phi_i)_{\Omega}$ \\

$M_h$ & $N\times N$ & $\left(\bu_h,\bv_h\right)$ & $(M_h)_{ij} = (\boldsymbol\phi_j,\boldsymbol\phi_i)_{\Omega}$ \\

$M_B$ & $N\times N_B$ & $\left\langle\bg,(\bn\times\bv_h)\times\bn\right\rangle$ & $(M_B)_{ij} = \langle \boldsymbol\psi_j, (\boldsymbol\phi_i)_T \rangle_{\partial\Omega}$ \\

\bottomrule
  \end{tabular}
\end{center}
\end{table}
For $\eta\in\Gamma$, define the two discrete solution operators
\begin{align}\label{def:R1hR0h}
R_{1,h} =\left(A_h^1-iM_h^B-\eta M_h^N\right)^{-1}M_B \,\,\text{and}\,\, R_{0,h} =\left(A_h^0-iM_h^B-\eta M_h\right)^{-1}M_B.
\end{align}
Hence, the discrete operator $T_h(\eta)$ is given by 
$$T_h(\eta):=M_B^t(R_{1,h}(\eta)-R_{0,h}(\eta))$$
and a discrete eigenpair $(\eta,\bg_h)$ is such that $T_h(\eta)\bg_h=0$ with $\bg_h\neq 0$.



For each retained contour domain, Beyn's method~\cite{beyn2012integral} is applied to compute the eigenvalues. It extracts the eigenvalues inside the contour by constructing a low-rank approximation of the spectral projection through contour integrals. Assume there exist $m\leq N_B$ eigenvalues $\{\lambda_j \}_{j=1}^m$ inside $\mathcal D$ and no eigenvalue lies on $\Gamma$. Let $V\in\mathbb{C}^{N_B\times \ell}$, $\ell\leq N_B$, be a random full rank matrix with $\ell \geq m$. Define 
\begin{align*}
C_0=\frac{1}{2\pi i}\int_{\Gamma}T_h(\eta)^{-1}V\,\mathrm{d}\eta\in\mathbb{C}^{N_B\times \ell}, \quad \text{and}\quad C_1=\frac{1}{2\pi i}\int_{\Gamma}\eta T_h(\eta)^{-1}V\,\mathrm{d}\eta\in\mathbb{C}^{N_B\times \ell}.
\end{align*}
The singular value decomposition of $C_0$ is used to determine the numerical rank, and the resulting reduced eigenvalue problem constructed from $C_0$ and $C_1$ yields the approximations of the eigenvalues and eigenvectors.

The computational procedure is summarized in Algorithm~\ref{alg:pmCIM}.

\begin{algorithm}\caption{pmCIM}\label{alg:pmCIM}
\begin{algorithmic}
\STATE{Input: domain $\O$; search region $\mathcal D$; numbers of quadrature points $N_{\rm SIM}$ and $N_{\rm Beyn}$; thresholds $tol_{\rm ind}$, $tol_{\rm svd}$ and $tol_{\rm res}$; number of probing vectors $\ell$.}
\begin{enumerate}
\item[1.]  Generate a tetrahedral mesh of $\O$ and assemble the matrices $A_h^1$, $M_h^B$, $M_h^N$, $A_h^0$, $M_h$, and $M_B$.
\item[2.] Partition $\mathcal D$ into squares $\{\mathcal{Q}_k\}_{k=1}^{I}$, and let $\{\mathcal{D}_k\}_{k=1}^{I}$ be their circumscribed disks with center $c_k$ and radius $r_k$ with $k=1,\cdots,I$.
\item[3.] \textbf{SIM screening.}
\begin{enumerate}
\item Choose a random $\bof_h\in V_h^B$;
\item For each disk $\mathcal D_k$, use the trapezoidal rule to approximate~\eqref{eq:proj} with $N_{SIM}$ nodes, choose $\eta_j=c+re^{i\theta_j}$, $\theta_j=2\pi(j-1)/N_{SIM},\, j=1,\cdots, N_{SIM}$ and compute the spectral indicator for $\mathcal{D}_k$ as follows:
\begin{align}\label{def:IndicatorID}
I_{\mathcal{D}_k} :=\left|\frac{1}{2\pi i}\sum_{j=1}^{N_{SIM}}\omega_j\,\left\langle\bof_h,T_h^{-1}(\eta_j)\bof_h\right\rangle\right|\quad \text{with} \quad \omega_j= \frac{i r_k e^{i\theta_j}2\pi}{N_{SIM}}.
\end{align}
\item Retain $\mathcal D_k$ if $I_{\mathcal D_k}>tol_{ind}$; otherwise discard it.
\end{enumerate}
\item[4.] \textbf{Beyn extraction on retained disks.}
\begin{enumerate}
\item Generate a random matrix $V\in \mathbb{C}^{N_B\times \ell}$. Let $\theta_j=2\pi(j-1)/N_{Beyn}$, $\eta_j=c+re^{i\theta_j}$ with $j=1,\cdots,N_{Beyn}$ and compute $C_{0,N_{Beyn}}$ and $C_{1,N_{Beyn}}$ by:
\begin{align}
C_0 &\approx C_{0,N_{Beyn}}=\frac{1}{N_{Beyn}}\sum_{j=1}^{N_{Beyn}}re^{i\theta_j}T_h(\eta_j)^{-1} V, \label{def:C0}\\
C_1 &\approx C_{1,N_{Beyn}}=\frac{1}{N_{Beyn}}\sum_{j=1}^{N_{Beyn}}re^{i\theta_j}\eta_jT_h(\eta_j)^{-1} V.\label{def:C1} 
\end{align}
\item Compute the SVD of $C_{0,N_{Beyn}}$, i.e., $C_{0,N_{Beyn}}=U\Sigma W^H$, where $U\in\mathbb C^{N_B\times \ell}$, $\Sigma=diag(\sigma_1,\sigma_2,\cdots,\sigma_{\ell})$ with $\sigma_1\geq\sigma_2\geq\cdots\geq\sigma_\ell\geq 0$, $W\in\mathbb C^{\ell\times\ell}$;
\item Determine the numerical rank $p$, $0<p<\ell$ by $\sigma_p>tol_{svd}\geq\sigma_{p+1}$. 
Set $U_p=U(:,1:p)$, $W_p=W(:,1:p)$ and $\Sigma_p=\operatorname{diag}(\sigma_1,\ldots,\sigma_p)$;
\item Form $G_i$ by $G_i:=U_p^HC_{1,N_{Beyn}}W_p\Sigma_p^{-1}\in\mathbb{C}^{p\times p}$ and compute the eigenvalues and eigenvectors of $G_i$, i.e., $G_i\boldsymbol s_j=\eta_j\boldsymbol s_j$. Then obtain the boundary eigenvectors of $T_h(\cdot)$ by $\boldsymbol g_{j,h}=U_p\boldsymbol s_j$;
\item Discard candidate eigenvalues that lie outside the square associated with $\mathcal D_i$.
\end{enumerate}

\item[5.] \textbf{Verification and compute eigenfunctions.}
\begin{enumerate}
\item For each candidate pair $(\eta_j,\boldsymbol g_{j,h})$, compute relative residual by ${\rm Res}(\eta_j,\bg_{j,h})=\|T_h(\eta_j)\bg_{j,h}\|_2/
\|\bg_{j,h}\|_2$
and accept the candidate if ${\rm Res}(\eta_j,\boldsymbol g_{j,h})<tol_{res}$;
\item For each final eigenpair $(\eta_j,\boldsymbol g_{j,h})$, reconstruct the discrete transmission eigenfunctions by $\bu_{j,h}=R_{1,h}(\eta_j)\bg_{j,h}$ and $\bv_{j,h}=R_{0,h}(\eta_j)\bg_{j,h}$.
\end{enumerate}

\end{enumerate}

\STATE{Output: wave numbers $k_{j,h}=\sqrt{\eta_j}$ and reconstructed fields $(\boldsymbol u_{j,h},\boldsymbol v_{j,h})$.}
\end{algorithmic}
\end{algorithm}

The computations are parallel at two levels. First, different contour domains are independent and can be processed concurrently. Second, for each contour, the linear systems associated with distinct quadrature nodes are independent.

\subsection{Numerical Experiments}\label{sec:NumExp}

We report some numerical results to assess the performance of the method. 

Let $\Omega_1=[0,1]^3$ and $\Omega_2=\left([-1,1]\times [-1,1]\times [0,1/2]\right)\backslash \left([0,1]\times[0,1]\times[0,1/2]\right)$.

\noindent Four sets of $\boldsymbol{A}$ and $\boldsymbol{N}$ are considered
\begin{align*}
&(1)\,\boldsymbol{A}(x)=2I,\,\,\boldsymbol{N}(x)=16I, \qquad\quad\qquad (2)\, \boldsymbol{A}(x)=F_2,\,\,\boldsymbol{N}(x)=F_1,\quad \\
&(3)\,\boldsymbol{A}(x)=F_2,\,\, \boldsymbol{N}(x,k)=(1+\delta k^2)F_1, \,\,\, (4)\,\boldsymbol{A}(x)=2I,\,\, \boldsymbol{N}(x,k)=(16+\sin(kx_1))I,
\end{align*}
where
\begin{align*}
\small
F_1 =\begin{bmatrix}
16 & x_1 & x_2 \\
x_1 & 16 & x_3 \\
x_2 & x_3 & 14
\end{bmatrix},\,\, F_2 = \begin{bmatrix}
-\frac{x_1^2}{8}+\frac{9x_1}{8}-\frac{9x_2}{8}+\frac{65}{4} & \frac{3x_1}{8}-\frac{3x_2}{8}-\frac{3x_1^2}{8}+\frac{3}{4} & 0\\
\frac{3x_1^2}{8}-\frac{3x_1}{8}+\frac{3x_2}{8}-\frac{3}{4} & \frac{9x_1^2}{8}-\frac{x_1}{8}+\frac{x_2}{8}+\frac{55}{4} & 0 \\
0 & 0 & x_3^2+12
\end{bmatrix}.
\end{align*}
The settings of $\bA$ and $\bN$ satisfies the requirements in Section \ref{sec:int}. We note that the last two examples are dispersive models.

We generate a coarse tetrahedral mesh and uniformly refine it to obtain a sequence of meshes. The lowest-order N\'ed\'elec edge element is used for the discretization of the $H(\text{curl})$ space. Let $k_{h_i}$, $i=1,\cdots,L$ be the computed sequence approximating an eigenvalue $k$. The relative error and the convergence order of the eigenvalues are defined by
\begin{align*}
\text{err}_k^i=\frac{|k_{h_i}-k_{h_{i-1}}|}{|k_{h_{i-1}}|},\,\,i=2,\cdots,L \quad \text{and}\quad\text{order}_k^i=\frac{\log(\text{err}_k^i/\text{err}_k^{i-1})}{\log(h_i/h_{i-1})},\,\, i=3,\cdots,L.
\end{align*}

Let $(\bu_{h_i},\bv_{h_i})$ be the computed eigenfunction pair associated with
the eigenvalue $k$ on the mesh with size $h_i$, $i=1,\ldots,L$. We define the relative $L^2$ error and its corresponding convergence order as follows
\begin{align*}
{\rm err}^{uv,i}_{L^2} &= \frac{ \left( \|\bu_{h_i}-\bu_{h_{i-1}}\|_{L^2(\Omega)}^2 + \|\bv_{h_i}-\bv_{h_{i-1}}\|_{L^2(\Omega)}^2 \right)^{1/2} }{ \left( \|\bu_{h_{i-1}}\|_{L^2(\Omega)}^2 + \|\bv_{h_{i-1}}\|_{L^2(\Omega)}^2 \right)^{1/2}}, \qquad i=2,\ldots,L,\\
{\rm order}^{uv,i}_{L^2} &= \frac{ \log\left({\rm err}^{uv,i}_{L^2}/{\rm err}^{uv,i-1}_{L^2}\right) }{ \log(h_i/h_{i-1}) }, \hspace{2.0in} i=3,\ldots,L .
\end{align*}
Similarly, the relative $H({\rm curl})$ error is computed by replacing the $L^2$ norm in the above definition with the $H({\rm curl})$ norm, where $\|\bw\|_{H({\rm curl},\Omega)}^2 = \|\bw\|_{L^2(\Omega)}^2 + \|\nabla\times \bw\|_{L^2(\Omega)}^2$.

During computation, all searches are carried out in the $\eta$-plane, where $\eta=k^2$. For each search rectangle, a uniform subdivision is first generated, and the circumscribing disks of the subregions are used as contour domains. In all computations, the same basic pmCIM parameters are used, with only minor case-dependent adjustments when necessary. The search rectangle is divided into $4\times 4$ subregions. The spectral indicator step uses $tol_{\rm ind}=10^{-6}$ and $N_{\rm SIM}=64$ quadrature points. In Beyn's method, we take $N_{\rm Beyn}=30$ and $\ell=3$. The singular value threshold is set to ${tol}_{\rm svd}=10^{-10}$, and the residual tolerance in the verification step is ${tol}_{\rm res}=10^{-2}$. The search regions used for different domains and coefficients are listed in Table~\ref{tab:pmcim_settings}.

\begin{table}[!htbp]
\centering
\caption{Parameters settings for different domains and coefficients.}
\small
\renewcommand{\arraystretch}{1.10}
\setlength{\tabcolsep}{12pt}

\begin{tabular}{c c c}
\toprule
Domain & Case & Search region\\
\midrule

\multirow{3}{*}{$\Omega_1$}
& $\bA=2I,\ \bN=16I$ & $[1.2,\,2.4]\times[-0.6,\,0.6]$ \\

& $\bA=F_2,\ \bN=F_1$ & $[12,\,20]\times[-4,\,4]$ \\

& $\bA=F_2,\ \bN=(1+\delta k^2)F_1$ & $[10,\,22]\times[-6,\,6]$ \\

& $\bA=2I,\ \bN=(16+\sin(kx_1))I$ & $[1,\,3]\times[-1,\,1]$ \\

\midrule

\multirow{3}{*}{$\Omega_2$}
& $\bA=2I,\ \bN=16I$ & $[0.3,\,1.3]\times[-0.3,\,0.3]$ \\

& $\bA=F_2,\ \bN=F_1$  & $[6,\,12]\times[-3,\,3]$ \\

& $\bA=F_2,\ \bN=(1+\delta k^2)F_1$ & $[3,\,13]\times[-5,\,5]$ \\

& $\bA=2I,\ \bN=(16+\sin(kx_1))I$ & $[0.1,\,2.1]\times[-1,\,1]$ \\

\bottomrule
\end{tabular}
\label{tab:pmcim_settings}
\end{table}

\begin{table}[!t]
\centering
\caption{Convergence results for $\Omega_1$ under different coefficients.}
%
\footnotesize
\renewcommand{\arraystretch}{0.85}
\setlength{\tabcolsep}{4pt}

\resizebox{\textwidth}{!}{
\begin{tabular}{c c c c c c c c c}
\toprule
Eigenpair
& $h$ & $k_h$ & ${\rm err}_k$ & ${\rm order}_k$
& ${\rm err}^{uv}_{L^2}$ & ${\rm order}_{L^2}$
& ${\rm err}^{uv}_{H({\rm curl})}$ & ${\rm order}_{H({\rm curl})}$ \\
\midrule

\multicolumn{9}{c}{Case 1: $\bA=2I,\ \bN=16I$} \\
\midrule

\multirow{5}{*}{$\#1$}
& $1$ & $1.2625$ & -- & -- & -- & -- & -- & -- \\
& $1/2$ & $1.2260$ & $2.89\times 10^{-2}$ & -- & $9.90\times 10^{-1}$ & -- & $9.91\times 10^{-1}$ & -- \\
& $1/4$ & $1.2079$ & $1.48\times 10^{-2}$ & $0.966$ & $4.96\times 10^{-1}$ & $0.997$ & $4.88\times 10^{-1}$ & $1.022$ \\
& $1/8$ & $1.2085$ & $5.65\times 10^{-4}$ & $4.712$ & $4.01\times 10^{-1}$ & $0.307$ & $3.89\times 10^{-1}$ & $0.327$ \\

\midrule

\multirow{5}{*}{$\#2$}
& $1$ & $1.5250$ & -- & -- & -- & -- & -- & -- \\
& $1/2$ & $1.3969$ & $8.40\times 10^{-2}$ & -- & $8.46\times 10^{-1}$ & -- & $8.80\times 10^{-1}$ & -- \\
& $1/4$ & $1.4514$ & $3.90\times 10^{-2}$ & $1.106$ & $9.17\times 10^{-1}$ & $-0.117$ & $9.19\times 10^{-1}$ & $-0.063$ \\
& $1/8$ & $1.4660$ & $1.00\times 10^{-2}$ & $1.960$ & $8.21\times 10^{-1}$ & $0.161$ & $8.21\times 10^{-1}$ & $0.163$ \\

\midrule

\multicolumn{9}{c}{Case 2: $\bA=F_2,\ \bN=F_1$} \\
\midrule

\multirow{4}{*}{$\#1$}
& $1/2$ & $4.0045$ & -- & -- & -- & -- & -- & -- \\
& $1/4$ & $3.9100$ & $2.36\times 10^{-2}$ & -- & $9.99\times 10^{-1}$ & -- & $9.98\times 10^{-1}$ & -- \\
& $1/8$ & $3.8760$ & $8.69\times 10^{-3}$ & $1.441$ & $1.99\times 10^{-1}$ & $2.329$ & $2.12\times 10^{-1}$ & $2.237$ \\
& $1/16$ & $3.8675$ & $2.18\times 10^{-3}$ & $1.993$ & $1.01\times 10^{-1}$ & $0.978$ & $1.08\times 10^{-1}$ & $0.975$ \\

\midrule

\multirow{4}{*}{$\#2$}
& $1/2$ & $4.1311$ & -- & -- & -- & -- & -- & -- \\
& $1/4$ & $4.3973$ & $6.44\times 10^{-2}$ & -- & $9.85\times 10^{-1}$ & -- & $9.76\times 10^{-1}$ & -- \\
& $1/8$ & $4.4559$ & $1.33\times 10^{-2}$ & $2.272$ & $9.90\times 10^{-1}$ & $-0.007$ & $9.93\times 10^{-1}$ & $-0.025$ \\
& $1/16$ & $4.4248$ & $6.99\times 10^{-3}$ & $0.933$ & $5.11\times 10^{-1}$ & $0.954$ & $4.85\times 10^{-1}$ & $1.034$ \\

\midrule
  
\multicolumn{9}{c}{Case 3: $\bA=F_2,\ \bN=(1+\delta k^2)F_1$ with $\delta =0.1$} \\
\midrule

\multirow{4}{*}{$\#1$}
& $1/2$ & $3.9954$ & -- & -- & -- & -- & -- & -- \\
& $1/4$ & $3.7582$ & $6.05\times 10^{-2}$ & -- & $9.99\times 10^{-1}$ & -- & $9.99\times 10^{-1}$ & -- \\
& $1/8$ & $3.5084$ & $8.32\times 10^{-2}$ & $-0.460$ & $1.00\times 10^{0}$ & $-0.002$ & $1.00\times 10^{0}$ & $-0.001$ \\
& $1/16$ & $3.5312$ & $4.45\times 10^{-2}$ & $0.903$ & $9.94\times 10^{-1}$ & $0.009$ & $9.94\times 10^{-1}$ & $0.008$ \\

\midrule

\multirow{4}{*}{$\#2$}
& $1/2$ & $4.1827$ & -- & -- & -- & -- & -- & -- \\
& $1/4$ & $4.5779$ & $1.20\times 10^{-1}$ & -- & $9.99\times 10^{-1}$ & -- & $9.99\times 10^{-1}$ & -- \\
& $1/8$ & $4.3515$ & $6.93\times 10^{-2}$ & $0.795$ & $1.00\times 10^{0}$ & $-0.001$ & $9.99\times 10^{-1}$ & $-0.001$ \\
& $1/16$ & $4.4683$ & $1.00\times 10^{0}$ & $1.328$ & $9.86\times 10^{-1}$ & $0.020$ & $9.93\times 10^{-1}$ & $0.010$ \\

\midrule

\multicolumn{9}{c}{Case 4: $\bA=2I,\ \bN=(16+\sin(kx_1))I$} \\
\midrule

\multirow{4}{*}{$\#1$}
& $1/2$ & $1.0451$ & -- & -- & -- & -- & -- & -- \\
& $1/4$ & $1.0237$ & $2.05\times 10^{-2}$ & -- & $9.99\times 10^{-1}$ & -- & $9.99\times 10^{-1}$ & -- \\
& $1/8$ & $1.0701$ & $4.53\times 10^{-2}$ & $-1.146$ & $9.75\times 10^{-1}$ & $0.035$ & $9.74\times 10^{-1}$ & $0.037$ \\
& $1/16$ & $1.0793$ & $8.65\times 10^{-3}$ & $2.390$ & $9.61\times 10^{-1}$ & $0.021$ & $9.54\times 10^{-1}$ & $0.030$ \\

\midrule

\multirow{4}{*}{$\#2$}
& $1/2$ & $1.2903$ & -- & -- & -- & -- & -- & -- \\
& $1/4$ & $1.2445$ & $3.55\times 10^{-2}$ & -- & $9.98\times 10^{-1}$ & -- & $9.98\times 10^{-1}$ & -- \\
& $1/8$ & $1.3361$ & $7.36\times 10^{-2}$ & $-1.054$ & $9.86\times 10^{-1}$ & $0.018$ & $9.85\times 10^{-1}$ & $0.019$ \\
& $1/16$ & $1.2772$ & $4.41\times 10^{-2}$ & $0.740$ & $8.40\times 10^{-1}$ & $0.231$ & $8.61\times 10^{-1}$ & $0.194$ \\

\bottomrule
\end{tabular}
}
\label{tab:domain1_eigenpairs}
\end{table}

The numerical results for $\Omega_1$ are listed in Table~\ref{tab:domain1_eigenpairs}. Even for the coarse meshes considered here, the resulting finite element systems are already large. We reported the first two eigenpairs under the four coefficient settings. Overall, the computed eigenvalues show a clear stabilization as the mesh is refined. Although the convergence rates are not always regular on the coarser meshes, they gradually approach a stable behavior as the mesh size decreases. 
With further mesh refinement, both eigenvalues and eigenfunctions show consistent convergence behavior.

The results for the nonconvex domain $\Omega_2$ are presented in Table~\ref{tab:domain2_eigenpairs}. Compared with the unit cube, the convergence is generally slower and less regular, particularly for the eigenfunctions. This behavior is consistent with the reduced regularity of Maxwell fields on nonconvex polyhedral domains. Despite this additional difficulty, the computed eigenvalues become progressively more stable under mesh refinement for all four coefficient settings. 

\begin{table}[!t]
\centering
\caption{Convergence results for $\Omega_2$ under different coefficients.}
\resizebox{\textwidth}{!}{
\begin{tabular}{c c c c c c c c c}
\toprule
Eigenpair
& $h$ & $k_h$ & ${\rm err}_k$ & ${\rm order}_k$
& ${\rm err}^{uv}_{L^2}$ & ${\rm order}_{L^2}$
& ${\rm err}^{uv}_{H({\rm curl})}$ & ${\rm order}_{H({\rm curl})}$ \\
\midrule

\multicolumn{9}{c}{Case 1: $\bA=2I,\ \bN=16I$} \\
\midrule

\multirow{5}{*}{$\#1$}
& $1$ & $0.8674$ & -- & -- & -- & -- & -- & -- \\
& $1/2$ & $0.7653$ & $1.18\times 10^{-1}$ & -- & $1.00\times 10^{0}$ & -- & $1.00\times 10^{0}$ & -- \\
& $1/4$ & $0.7962$ & $4.03\times 10^{-2}$ & $1.545$ & $3.73\times 10^{-1}$ & $1.422$ & $4.88\times 10^{-1}$ & $1.849$ \\
& $1/8$ & $0.8100$ & $1.74\times 10^{-2}$ & $1.213$ & $2.32\times 10^{-1}$ & $0.689$ & $7.23\times 10^{-1}$ & $0.858$ \\

\midrule

\multirow{5}{*}{$\#2$}
& $1$ & $1.0250$ & -- & -- & -- & -- & -- & -- \\
& $1/2$ & $1.0364$ & $1.12\times 10^{-2}$ & -- & $5.11\times 10^{-1}$ & -- & $5.14\times 10^{-1}$ & -- \\
& $1/4$ & $1.0381$ & $1.61\times 10^{-3}$ & $2.798$ & $4.14\times 10^{-1}$ & $0.304$ & $4.09\times 10^{-1}$ & $0.330$ \\
& $1/8$ & $1.0388$ & $6.74\times 10^{-4}$ & $1.256$ & $2.15\times 10^{-1}$ & $0.941$ & $2.01\times 10^{-1}$ & $1.026$ \\

\midrule

\multicolumn{9}{c}{Case 2: $\bA=F_2,\ \bN=F_1$} \\
\midrule

\multirow{4}{*}{$\#1$}
&  $1$  & $3.2486$ & -- & -- & -- & -- & -- & -- \\
& $1/2$ & $2.5355$ & $2.20\times 10^{-1}$ & -- & $1.00\times 10^{0}$ & -- & $1.00\times 10^{0}$ & -- \\
& $1/4$ & $2.6172$ & $3.22\times 10^{-2}$ & $2.768$ & $3.75\times 10^{-1}$ & $1.415$ & $9.98\times 10^{-1}$ & $1.788$ \\
& $1/8$ & $2.6527$ & $1.36\times 10^{-2}$ & $1.248$ & $2.28\times 10^{-1}$ & $0.718$ & $2.12\times 10^{-1}$ & $0.859$ \\

\midrule

\multirow{4}{*}{$\#2$}
&  $1$  & $3.4556$ & -- & -- & -- & -- & -- & -- \\
& $1/2$ & $3.2339$ & $6.41\times 10^{-2}$ & -- & $1.00\times 10^{0}$ & -- & $1.00\times 10^{0}$ & -- \\
& $1/4$ & $3.1456$ & $2.73\times 10^{-2}$ & $1.232$ & $3.30\times 10^{-1}$ & $1.601$ & $3.91\times 10^{-1}$ & $1.355$ \\
& $1/8$ & $3.1150$ & $9.73\times 10^{-3}$ & $1.488$ & $1.58\times 10^{-1}$ & $1.060$ & $2.05\times 10^{-1}$ & $0.932$ \\

\midrule

\multicolumn{9}{c}{Case 3: $\bA=F_2,\ \bN=(1+\delta k^2)F_1$ with $\delta=0.1$} \\
\midrule

\multirow{4}{*}{$\#1$}
&  $1$  & $1.9387$ & -- & -- & -- & -- & -- & -- \\
& $1/2$ & $2.0935$ & $7.99\times 10^{-2}$ & -- & $5.07\times 10^{-1}$ & -- & $3.35\times 10^{-1}$ & -- \\
& $1/4$ & $2.1447$ & $2.44\times 10^{-2}$ & $1.708$ & $3.76\times 10^{-1}$ & $0.434$ & $2.86\times 10^{-1}$ & $0.227$ \\
& $1/8$ & $2.1667$ & $1.03\times 10^{-2}$ & $1.252$ & $2.29\times 10^{-1}$ & $0.712$ & $1.57\times 10^{-1}$ & $0.863$ \\

\midrule

\multirow{4}{*}{$\#2$}
&  $1$  & $2.5113$ & -- & -- & -- & -- & -- & -- \\
& $1/2$ & $2.9673$ & $1.82\times 10^{-1}$ & -- & $1.00\times 10^{0}$ & -- & $1.00\times 10^{0}$ & -- \\
& $1/4$ & $2.5755$ & $1.32\times 10^{-1}$ & $0.460$ & $1.00\times 10^{0}$ & $0.000$ & $1.00\times 10^{-1}$ & $0.000$ \\
& $1/8$ & $2.5718$ & $1.46\times 10^{-3}$ & $6.502$ & $2.14\times 10^{-1}$ & $2.222$ & $2.06\times 10^{-1}$ & $2.279$ \\

\midrule

\multicolumn{9}{c}{Case 4: $\bA=2I,\ \bN=(16+\sin(kx_1))I$} \\
\midrule

\multirow{4}{*}{$\#1$}
& $1$ & $0.8671$ & -- & -- & -- & -- & -- & -- \\
& $1/2$ & $0.7701$ & $1.12\times 10^{-1}$ & -- & $1.00\times 10^{0}$ & -- & $1.00\times 10^{0}$ & -- \\
& $1/4$ & $0.8015$ & $4.07\times 10^{-2}$ & $1.459$ & $3.72\times 10^{-1}$ & $1.425$ & $2.77\times 10^{-1}$ & $1.851$ \\
& $1/8$ & $0.8156$ & $1.76\times 10^{-2}$ & $1.212$ & $2.31\times 10^{-1}$ & $0.689$ & $1.53\times 10^{-1}$ & $0.858$ \\

\midrule

\multirow{4}{*}{$\#2$}
& $1$ & $1.0276$ & -- & -- & -- & -- & -- & -- \\
& $1/2$ & $1.0414$ & $1.34\times 10^{-2}$ & -- & $5.03\times 10^{-1}$ & -- & $5.00\times 10^{-1}$ & -- \\
& $1/4$ & $1.0385$ & $2.80\times 10^{-3}$ & $2.262$ & $4.10\times 10^{-1}$ & $0.296$ & $4.06\times 10^{0}$ & $0.301$ \\
& $1/8$ & $1.0419$ & $3.26\times 10^{-3}$ & $-0.217$ & $2.15\times 10^{-1}$ & $0.931$ & $2.01\times 10^{-1}$ & $1.018$ \\

\bottomrule
\end{tabular}
}
\label{tab:domain2_eigenpairs}
\end{table}

\begin{figure}[htbp]
\centering
\includegraphics[width=0.95\textwidth]{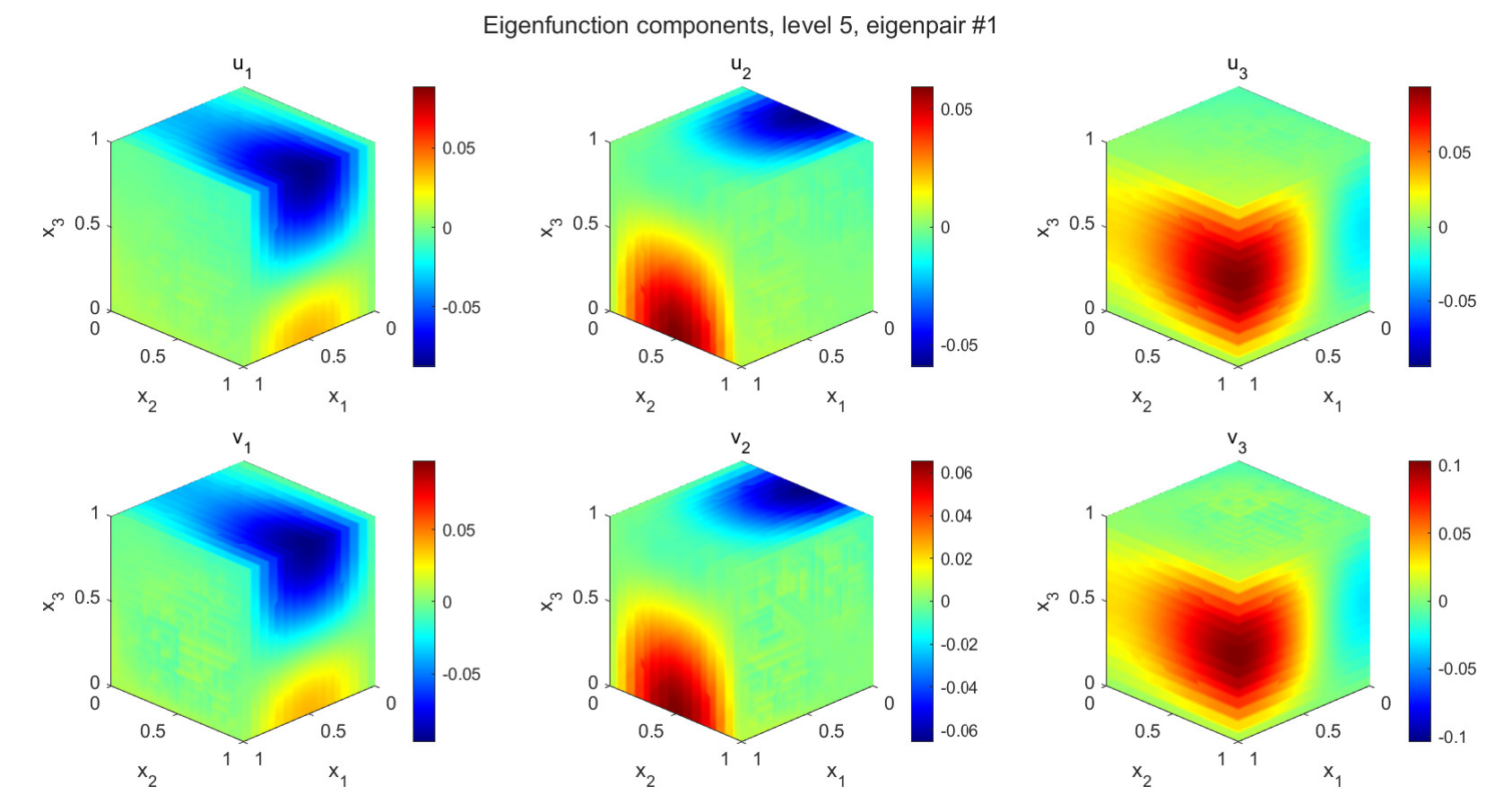}
\caption{Component-wise plots of the real parts of the first computed eigenfunction pairs $(\bu_h,\bv_h)$ for $\Omega_1$ with $A=2I$ and $N=16I$ on the finest mesh. 
}
\label{fig:O1_case1_eigfun}
\end{figure}

\begin{figure}[htbp]
\centering
\includegraphics[width=0.95\textwidth]{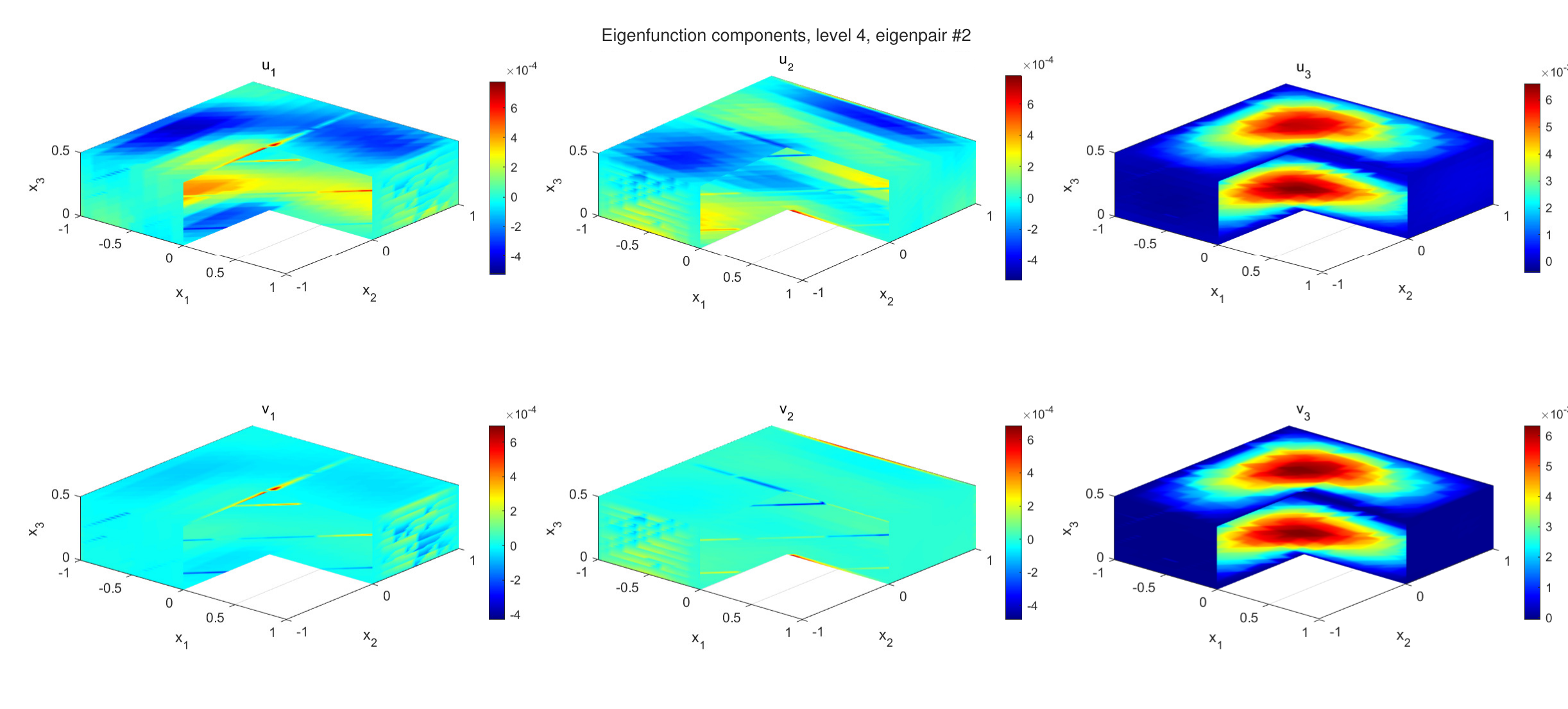}
\caption{Component-wise plots of the real parts of the second two computed eigenfunction pairs $(\bu_h,\bv_h)$ for $\Omega_2$ with $A=F_2$ and $N=(1+\delta k^2)F_1$ with $\delta=0.1$ on the finest mesh. 
}
\label{fig:O2_case2_eigfun}
\end{figure}

To further illustrate the computed eigenfunctions, we plot the real parts of the components of some recovered eigenfunction pairs on the finest mesh as shown in Figure~\ref{fig:O1_case1_eigfun} $\sim$ Figure \ref{fig:O2_case2_eigfun}. It can be observed that different components may have significantly different magnitudes, and some components are relatively small compared with the dominant ones. These small-amplitude components reflect the intrinsic modal structures of the corresponding transmission eigenfunctions. 

\begin{figure}[htbp]
\centering
\includegraphics[width=0.95\textwidth]{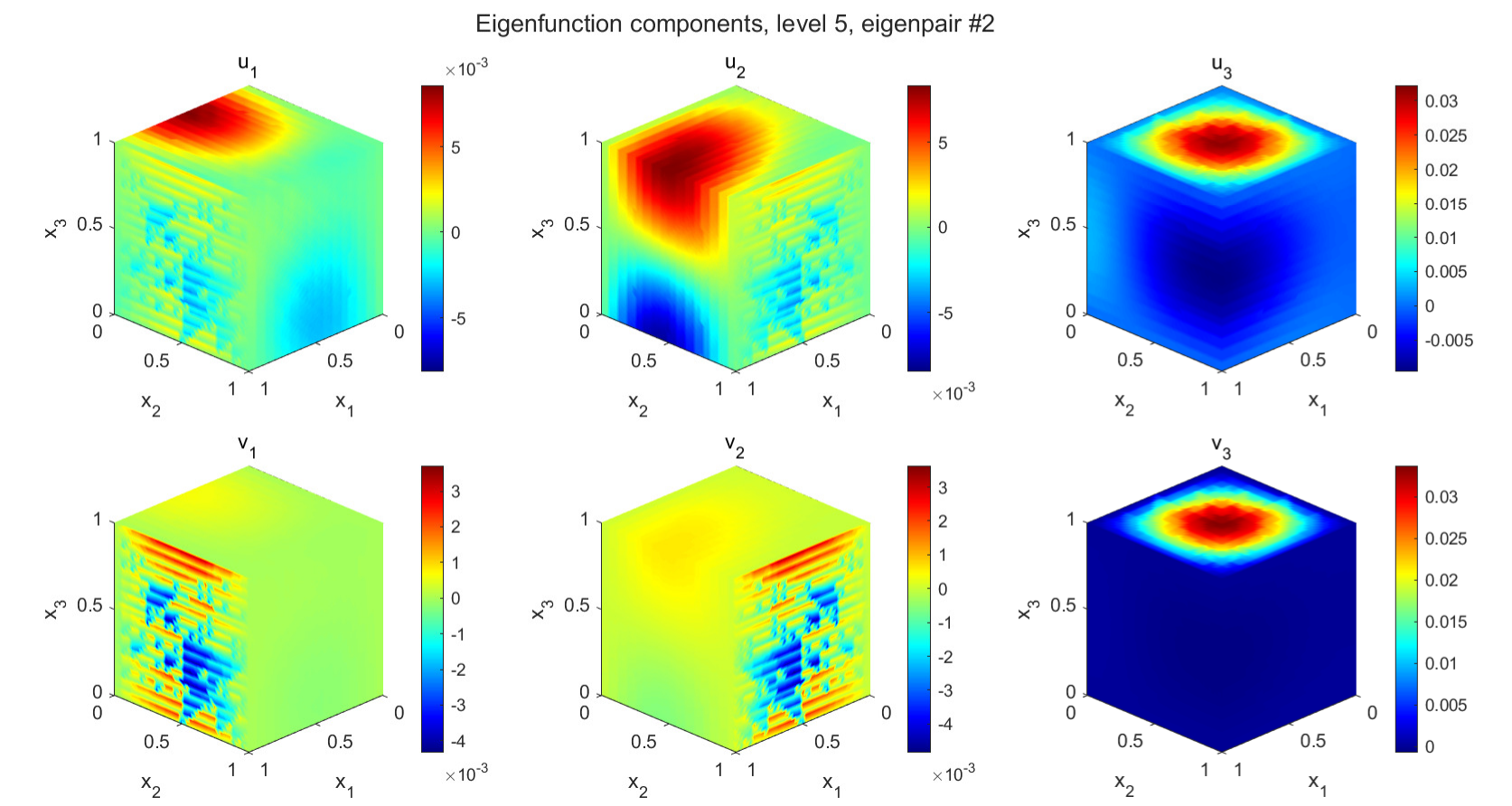}
\caption{Component-wise plots of the real parts of the second two computed eigenfunction pairs $(\bu_h,\bv_h)$ for $\Omega_1$ with $A=F_2$ and $N=F_1$ on the finest mesh. 
}
\label{fig:O1_case3_eigfun}
\end{figure}
%
%
%
\begin{figure}[htbp]
\centering
\includegraphics[width=1.0\textwidth]{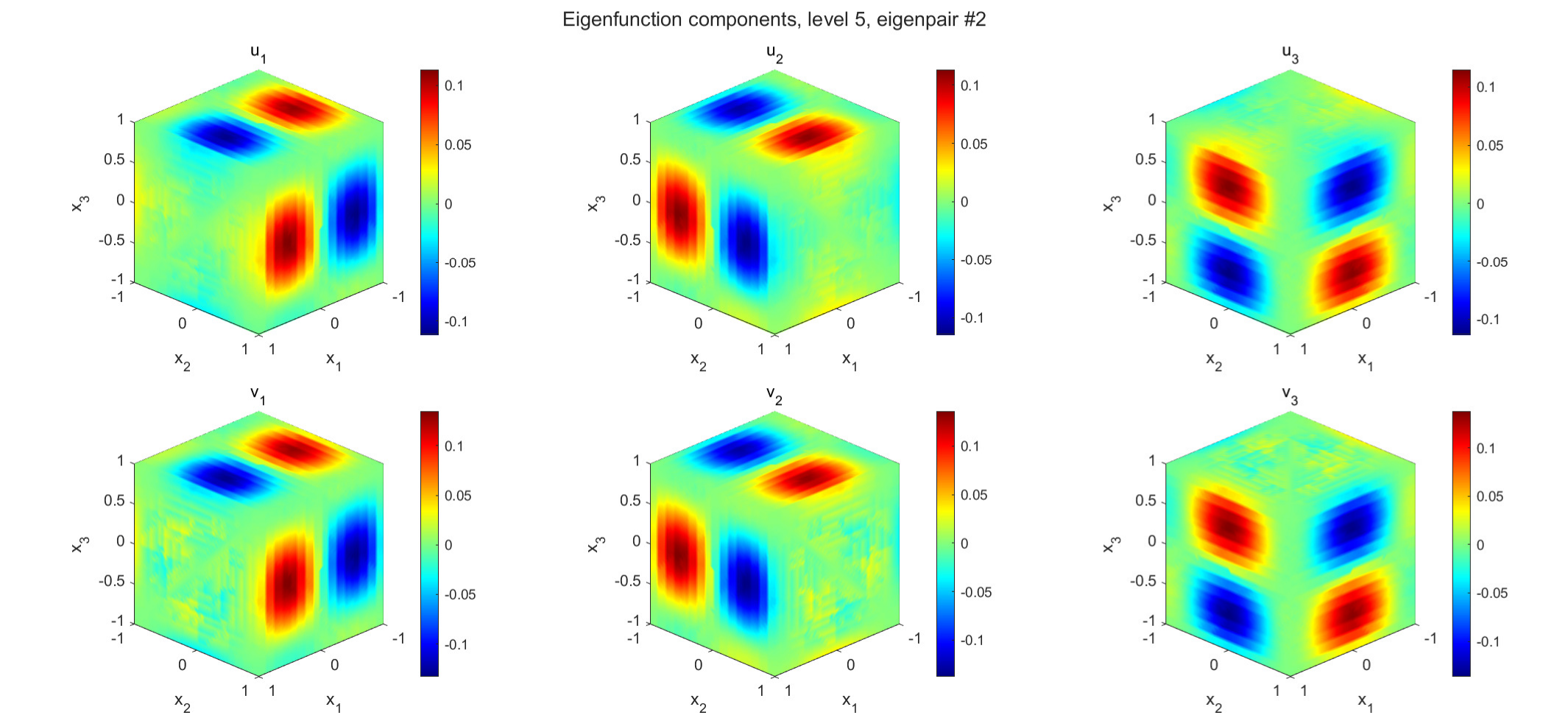}
\caption{Component-wise plots of the real parts of the second two computed eigenfunction pairs $(\bu_h,\bv_h)$ for $\Omega_1$ with $A=2I$ and $N=(16+\sin(kx_1))I$ on the finest mesh.
}
\label{fig:O1_case2_eigfun}
\end{figure}
%
%
%
%
%
%
%

\section{Conclusions}\label{sec:con} %

In this work, we study the Maxwell transmission eigenvalue problem in dispersive media by reformulating it as an eigenvalue problem for a holomorphic Fredholm operator function of index zero. Using Nédélec edge elements, we build a conforming finite element discretization and establish convergence of the approximated eigenvalues based on spectral approximation theory for holomorphic Fredholm operator functions. Moreover, a parallel multi-step contour integral method (pmCIM) is employed to compute eigenvalues in prescribed complex regions efficiently. Numerical results in three dimensions validate the accuracy and robustness of the proposed formulation and algorithm. Future work will focus on more efficient preconditioned solvers. In addition, we will investigate more general dispersive material models and establish the Fredholm property of the associated boundary operator functions beyond the assumptions considered in this work.

\section*{Acknowledgments}
The author would like to thank Prof. Sun Jiguang and Prof. Cui Jintao for helpful comments and suggestions on this work.

\bibliographystyle{siamplain}
\bibliography{Ref}
\end{document}